\documentclass[11pt]{article}

\usepackage[utf8]{inputenc}
\usepackage{graphicx}
\graphicspath{ {images/} }
\usepackage[a4paper,width=150mm,top=25mm,bottom=25mm]{geometry}

\usepackage{xcolor}

\usepackage{amsmath}
\usepackage{amssymb}
\usepackage{mathtools}

\usepackage{amsthm}

\usepackage{stackengine}
\stackMath

\usepackage{braket}
\usepackage{xfrac}

\usepackage{tikz-cd}

\usetikzlibrary{ shapes }

\usepackage{hyperref}
\usepackage[noabbrev]{cleveref}

\usepackage[font={small}]{caption}

\usepackage{verbatim}

\usepackage{mathrsfs}

\usepackage[cal=boondox,scr=boondoxo]{mathalfa}

\usepackage{enumerate}

\usepackage[numbers,sort]{natbib}

\usepackage{datetime}

\usepackage{blkarray}
\usepackage{nicematrix}
\NiceMatrixOptions{
code-for-first-row = \color{black} ,
code-for-last-row = \color{black} ,
code-for-first-col = \color{black} ,
code-for-last-col = \color{black}
}

\usetikzlibrary{patterns}

\usepackage[export]{adjustbox}
\usepackage{mwe}

\usepackage{authblk}

\usepackage{pagecolor}

\usepackage{standalone}
\pgfdeclarelayer{bg} 
\pgfdeclarelayer{backgroundlayer}
\pgfsetlayers{backgroundlayer,bg,main}
\usetikzlibrary{calc}
\newtheorem{numerable}{Numerable}
\newtheorem{lemma}[numerable]{Lemma}

\newtheorem{theorem}[numerable]{Theorem}
\newtheorem{corollary}[numerable]{Corollary}

\numberwithin{numerable}{section}
\theoremstyle{definition}
\newtheorem{definition}[numerable]{Definition}
\newtheorem{example}[numerable]{Example}

\definecolor{details}{RGB}{0,0,255}
\definecolor{task}{RGB}{0,191,0}
\definecolor{sketch}{RGB}{255,0,0}

\newcommand{\N}{\mathbb{N}}

\newcommand{\Z}{\mathbb{Z}}

\newcommand{\R}{\mathbb{R}}

\newcommand{\T}{\mathbb{T}}
\newcommand{\TP}{\mathbb{T}\mathbb{P}}

\newcommand{\conv}{\operatorname{conv}}

\newcommand{\im}{\operatorname{im}}
\newcommand{\pos}{\operatorname{pos}}
\newcommand{\vol}{\operatorname{vol}}
\newcommand{\mns}{\texttt{mns}}
\newcommand{\tns}{\texttt{tns}}

\renewcommand{\min}{\operatorname{min}}
\renewcommand{\l}{\ell}

\newdateformat{numberdate}{\THEDAY.\THEMONTH.\THEYEAR}

\makeatletter
\patchcmd{\@part}{\markboth{}{}}{\markboth{#1}{}}{\typeout{done patching part!}}{\typeout{oh dear! could not patch part command...}}
\patchcmd{\@spart}{\nobreak}{\nobreak\markboth{#1}{}}{\typeout{done patching starred part!}}{\typeout{oh dear! could not patch starred part command...}}
\makeatother

\title{
    The tropical and zonotopal geometry of periodic timetables
}

\author{Enrico Bortoletto\thanks{Corresponding author; \texttt{bortoletto@zib.de}}}
\author{Niels Lindner} 
\author{Berenike Masing}
\affil{Zuse Institute Berlin, Takustr. 7, 14195 Berlin, Germany}

\date{Dated: \today}

\begin{document}

\maketitle

\begin{abstract}
    
    
    The Periodic Event Scheduling Problem (PESP) is the standard mathematical tool for optimizing periodic timetabling problems in public transport.
    A solution to PESP consists of three parts: a periodic timetable, a periodic tension, and integer periodic offset values.
    While the space of periodic tension has received much attention in the past, we explore geometric properties of the other two components, establishing novel connections between periodic timetabling and discrete geometry.
    Firstly, we study the space of feasible periodic timetables, and decompose it into polytropes, i.e., polytopes that are convex both classically and in the sense of tropical geometry.
    We then study this decomposition and use it to outline a new heuristic for PESP, based on the tropical neighbourhood of the polytropes.
    Secondly, we recognize that the space of fractional cycle offsets is in fact a zonotope.
    We relate its zonotopal tilings back to the hyperrectangle of fractional periodic tensions and to the tropical neighbourhood of the periodic timetable space.
    To conclude we also use this new understanding to give tight lower bounds on the minimum width of an integral cycle basis.
\end{abstract}

\section{Introduction}
\label{sec:introduction}

The timetable is the heart of a public transit system. 
Many public transportation networks across the world are operated in a periodic manner. 
The creation and optimization of periodic timetables is therefore an ubiquitous and frequent planning task. 
The standard mathematical model for periodic timetabling in public transport is the Periodic Event Scheduling Problem (PESP) developed in \cite{SerafiniUkovichPESP}. 
PESP is a challenging problem in various respects: On the theoretical side, the problem of finding a feasible periodic timetable is NP-hard for a given period time $T \geq 3$ \cite{Odijk1994, nachtigall_cutting_1996} or when the underlying constraint graph is series-parallel \cite{lindner_analysis_2022}. 
In practice, none of the current 22 benchmarking instances of the library \texttt{PESPlib} \cite{PespLib} could up to today be solved to proven optimality. 
However, a plenty of heuristic algorithms are available, connecting PESP with a zoo of well-known problems and techniques in combinatorial optimization, e.g., simplex algorithms \cite{nachtigall_solving_2008, goerigk_improving_2013}, maximum cuts \cite{lindner_new_2019}, matchings \cite{patzold_matching_2016}, and Boolean satisfiability \cite{grosmann_solving_2012, matos_solving_2020}. 
The prevailing approach to solve PESP instances exactly is mixed-integer programming \cite{LiebchenBook, borndorfer_concurrent_2020}. 
For example, the subway network of Berlin has been optimized by solving such a mixed-integer programming model \cite{liebchen_first_2008}. 

There are two central notions in periodic timetabling in public transport: A periodic timetable associates, intuitively speaking, a periodically repeating departure or arrival time to every stop of a trip. 
A periodic tension collects the durations of all activities in a public transport network, such as, e.g., driving between to neighboring stops, or transferring at a stop \cite{liebchen_modeling_2007}. 
A timetable can easily be computed from a tension and vice versa. 
In the context of mixed-integer programming, the convex hull of feasible periodic tensions is the central geometric object of study, and several classes of cutting places have been deduced by studying this polyhedron \cite{Odijk1994, nachtigall_cutting_1996, borndorfer_separation_2020, lindnerIntegerVertices}. We will formally introduce PESP and the relevant notions from periodic timetabling in \Cref{sec:context}.

The first aim of our paper is to study the geometry of the space $\Pi$ of feasible periodic timetables. 
It turns out that this space is intimately related with tropical geometry. 
The space $\Pi$ is naturally embedded into a torus and decomposes into pairwise disjoint polytropes, i.e., polytopes that are also a tropical convex hull of finitely many points \cite{JoswigKulas, JoswigLoho}. 
Analyzing the neighborhood relations of those polytropes, we outline a new primal heuristic for PESP. 
We will discuss the link between periodic timetables and tropical geometry in depth in \Cref{sec:tropical}.

A straightforward technique to model a periodicity constraint such as $Ax \equiv b \bmod T$ in a mixed-integer program is to rewrite the constraint as $Ax - b = Tz$ for some integral $z$. 
Therefore, the common mixed-integer programs for PESP make use of such modulo parameters. 
In \Cref{sec:zonotopes}, we will consider the space $Z$ of modulo parameters for the linear programming relaxation of the so-called cycle-based mixed-integer programming formulation for PESP. 
The space $Z$ is a zonotope, and as for a graphical zonotope, the maximal tiles of any fine zonotopal tiling of $Z$ correspond to spanning trees. Moreover, there is a certain duality to the space $\Pi$ of periodic timetables: In any fine zonotopal tiling of $Z$, the tiles containing lattice points correspond with vertices of the polytropes in the decomposition of $\Pi$. 
Conversely, the tropical vertices of those polytropes can be used to construct a fine zonotopal tiling of $Z$. 
We finish the section by connecting the zonotope $Z$ with the minimum width of integral cycle bases, a notion that can be used to estimate the efficiency of the cycle-based MIP formulation \cite{cyclePESPLiebchen}. 
As a byproduct, we obtain that the number of spanning trees of a graph is at most the product of the lengths of the cycles in any integral cycle basis.

We conclude the paper by giving an outlook in \Cref{sec:outlook}.

\section{Periodic Timetabling in Public Transport}
\label{sec:context}

\subsection{The Periodic Event Scheduling Problem}

We start by presenting the well-known standard formulation of the Periodic Event Scheduling Problem (PESP), first introduced by Serafini and Ukovich in \cite{SerafiniUkovichPESP}.
\begin{definition}
    An instance of the Periodic Event Scheduling Problem (PESP) consists of a directed graph $G$, a period $T \in \N$, lower and upper bounds $\l,u \in \R^{A(G)}$, and arc weights $w \in \R^{A(G)}$.
    A \emph{periodic timetable} is a vector $\pi \in \R^{V(G)}$, and any vector $x \in \R^{A(G)}$ such that $x_{ij} \equiv_T \pi_j - \pi_i$ is called a \emph{periodic tension} associated to $\pi$.
    A periodic tension $x$ is said to be \emph{feasible} for the given PESP instance if $\l \leq x \leq u$.
    A periodic timetable $\pi$ is said to be \emph{feasible} for the given PESP instance if there exists a feasible periodic tension $x$ associated to $\pi$.
    Given an instance as above, the \emph{Periodic Event Scheduling Problem (PESP)} consists in finding a feasible periodic timetable $\pi$ and a feasible associated tension $x$ such that the weighted tension $w^\top x$ is minimised.
\end{definition}

\begin{example}\label{ex:running}
    To illustrate the above definition, consider the small instance as given by the graph in \Cref{fig:ean}, which will serve as a running example. 
    \begin{figure}[htbp]
        \begin{center}
            \begin{tikzpicture}[scale=2.5]
					\tikzstyle{p} = [line width=1.5, ->]
					\tikzstyle{u} = [line width=1.5]
					\tikzstyle{q} = [line width=1.5, ->]
					\tikzstyle{r} = [line width=1.5, ->]
					\tikzstyle{v} = [draw, circle, inner sep=1, minimum width=12, font=\footnotesize]
					\tikzstyle{w} = [v]
					\tikzstyle{s} = [draw, rectangle, inner sep=1, minimum width=12, minimum height=12, font=\footnotesize]
					\tikzstyle{t} = [midway, font=\footnotesize]
					\tikzstyle{ta} = [t, above]
					\tikzstyle{tb} = [t, below]
					\tikzstyle{tr} = [t, right]
					\tikzstyle{tl} = [t, left]
					\tikzstyle{l} = [font=\scriptsize]
					\tikzstyle{ch1} = [blue]
					\tikzstyle{ch2} = [orange]
					\node[w, label=left:{\small{\textcolor{blue}{$\pi_0 = 0$}}}] (A) at (0, 0) {$v_0$ };
					\node[v,label=right:{\small{\textcolor{blue}{$\pi_1 = 8$}}}] (B) at (1, 1.5) {$v_1$};
					\node[w,label=right:{\small{\textcolor{blue}{$\pi_2 = 2$}}}] (C) at (2, 0) { $v_2$};
					
					\draw[p] (A) -- node[ta, black, sloped] {$[3,12]$} node[tb, ch2, sloped] {$x_{01} = 8$} (B);
					\draw[p] (B) -- node[ta, black, sloped] {$[4,13]$} node[tb, ch2, sloped] {$x_{12} = 4$} (C);
					\draw[p] (A) -- node[ta, black, sloped] {$[2,10]$} node[tb, ch2, sloped] {$x_{02} = 2$} (C);	
\end{tikzpicture}
        \end{center}
        \caption{
            Exemplary graph $G$ with labels $[l_{ij}, u_{ij}]$ for the arc-bounds for period time $T = 10$. 
            A feasible timetable and corresponding periodic tension are marked in color.
            }
        \label{fig:ean}
    \end{figure}
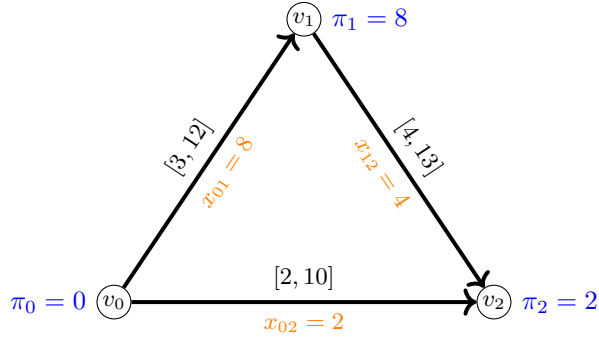
\end{example}

In practice $\l$ and $u$ are often integral. 
In this case, by a result of Odjik \cite{Odijk1994}, feasibility of the PESP instance implies the existence of an integral optimal periodic timetable and an associated integral periodic tension.
We will also make use of some standard assumptions \cite{LiebchenBook}: We suppose that $G$ is simple, $2$-connected, and has no pair of antiparallel arcs. 
Moreover, we assume $0 \leq \l < T$ and $0 \leq u - \l < T$.

\subsection{The Space of Feasible Periodic Timetables}

As a mixed-integer linear program the problem is formulated as 
\begin{equation}\label{eq:pesp-formulation}
    \begin{array}{ll@{}ll}
        \text{Minimise}     & w^\top x \\
        \text{subject to}   & \pi_j - \pi_i + T p_{ij} = x_{ij},& \quad  &\forall (i,j) \in A(G),\\
                            & \l \leq x \leq u, & \\
                            & p \in \Z^{A(G)},
    \end{array}
\end{equation}
where $p \in \Z^{A(G)}$ is called \emph{periodic offset}, a collection of auxiliary integer variable used to linearise the modulo constraints. 
Note that, with the incidence matrix $B$ of $G$, we can rewrite the constraints $\pi_j - \pi_i + T p_{ij} = x_{ij}$ for all $(i, j) \in A(G)$ as $-B^\top \pi + T p = x$.

Given that $\l \leq x$, the minimum weighted tension is bounded below by $w^\top \l$, which is why many authors prefer to translate the problem and consider the periodic slack $y \coloneqq x - \l$ instead.
In such cases we then have $y_{ij} \equiv_T \pi_j - \pi_i - \l_{ij}$ and $0 \leq y \leq u - \l$.

For example, expressing the feasible timetable from \Cref{fig:ean} in terms of slack, would result in $y = (y_{01}, y_{02}, y_{12}) = (5,0,0)$ with a periodic offset of $p = (p_{01},p_{02}, p_{12}) = (0,0, 1)$.

In view of the MIP \eqref{eq:pesp-formulation}, we can study the following two spaces:

\begin{definition}
    Given a PESP instance $(G, T, \l, u, w)$, we define the \emph{periodic tension polytope} as
    \begin{equation}
        X \coloneqq \conv\left\{ x \in \R^{A(G)} \middle| \
        \exists \pi \in \R^V \ \exists p \in \Z^{A(G)} \colon
        x = -B^\top \pi + T p, \
        \l \leq x \leq u
        \right\}.
    \end{equation}
\end{definition}

\begin{definition}
    Given a PESP instance $(G, T, \l, u, w)$, we define $\Pi$ as the space of feasible periodic timetables, i.e.,
    \begin{equation}
        \Pi \coloneqq \left\{ \pi \in \R^{V(G)} \middle|\ \exists p \in \Z^{A(G)}, \, \forall (i,j) \in A(G) \colon \l_{ij} \leq \pi_j - \pi_i + T p_{ij} \leq u_{ij}   \right\}.
    \end{equation}
\end{definition}

As per the introduction, the periodic tension polytope $X$ has been object of extensive study in periodic timetabling \cite{Odijk1994, nachtigall_cutting_1996, borndorfer_separation_2020, lindnerIntegerVertices}. 
However, the space $\Pi$ of feasible periodic timetables has not received much attention in the past. 
We will explore $\Pi$ with methods of tropical geometry in \Cref{sec:tropical}.

\subsection{The Space of Cycle Offsets}

Another well-known mixed-integer linear program, which we will make use of in \Cref{sec:zonotopes}, uses a cycle-based approach \cite{cyclePESPLiebchen}.

An \emph{oriented cycle} in $G$ is a vector $\gamma \in \{-1, 0, 1\}^{A(G)}$ such that $B\gamma = 0$, where $B$ is the incidence matrix of $G$. 
In other terms, $\gamma$ is a circuit in $G$, whose arcs can be traversed in forward or backward direction. 
The oriented cycles in $G$ generate a free abelian subgroup of $\Z^{A(G)}$, the \emph{cycle space} of $G$. The rank of the cycle space is the \emph{cyclomatic number} $\mu$, which for connected graphs can be computed as $\mu = |A(G)| - |V(G)| + 1$.

An \emph{integral cycle basis} $\mathscr B$ of $G$ is a $\mu$-tuple $(\gamma_1, \dots, \gamma_\mu)$ of oriented cycles in $G$ that generate the cycle space of $G$ as a $\Z$-module. 
The matrix $(\gamma_1, \dots, \gamma_\mu)^\top \in \Z^{\mathscr B \times A}$ is called the \emph{cycle matrix} of $\mathscr B$. 
A particular class of integral cycle bases is formed by fundamental cycle bases, which consist of the fundamental cycles of a spanning tree of $G$. 
We refer to \cite{kavitha_cycle_2009} for more background on cycle bases.

\begin{theorem}\label{thm:cycle-periodicity}
    Let $\mathscr B$ be an integral cycle basis of $G$ with cycle matrix $\Gamma$, and let $B$ denote the incidence matrix of $G$. 
    Then $\ker \Gamma = \im B^\top$ as $\Z$-modules. 
    In particular, a vector $x \in \mathbb R^{A(G)}$ is a periodic tension if and only if $\Gamma x \equiv_T 0$.
\end{theorem}

\begin{proof}
    The ``in particular'' statement is known as ``cycle periodicity property'' in the periodic timetabling literature \cite{cyclePESPLiebchen, nachtigallHabilitation} and follows from $\ker \Gamma = \im B^\top$ by taking the tensor product with $\Z/T\Z$. 
    The inclusion $\ker \Gamma \supseteq \im B^\top$ is clear, as the rows of $\Gamma$ belong to the kernel of $B$.
    Since $\Gamma$ has full row rank, the rank of $\ker \Gamma$ is $|A(G)| - \mu = |V(G)| - 1$, which coincides with the rank of the image of $B^\top$. 
    It follows that $\ker \Gamma / \im B^\top$ is a torsion group, and it remains to show that  $\ker \Gamma / \im B^\top = 0$.
    
    Let $x \in \ker \Gamma$. 
    Then there is a natural number $n$ such that $n \cdot x = B^\top \pi$ for some integral vector $\pi \in \Z^{V(G)}$. 
    Since the columns of $B^\top$ sum to zero, we can actually remove an arbitrary column of $B^\top$, so that $n \cdot x = B'^\top \pi'$ for some submatrix $B'$ of $B$ and some subvector $\pi'$ of $\pi$. 
    The matrix $B'^\top$ has now full column rank, so that we can choose an invertible submatrix $B''$ of $B$ such that $n \cdot x' = B''^\top \pi'$ for a suitable subvector $x'$ of $x$, and $\pi'$ is the unique solution to this system of linear equations. 
    By Cramer's rule, each entry of $\pi'$ is the quotient of a multiple of $n$ by the determinant of $B''$. 
    Since $B$ and hence $B''$ is totally unimodular, each entry of $\pi'$ is hence divisible by $n$. 
    Since this holds for any dropped column of $B^\top$, we conclude that $x = B^\top \frac{\pi}{n}$, and $\frac{\pi}{n} \in \Z^{V(G)}$.
\end{proof}

Let $\mathscr{B}$ be an integral cycle basis of the digraph $G$ and let $\Gamma$ be the corresponding cycle matrix. 
By \Cref{thm:cycle-periodicity}, PESP can be alternatively formulated as 
\begin{equation}\label{eq:pesp-formulation-cycle}
    \begin{array}{ll@{}ll}
        \text{Minimise}     & w^\top x \\
        \text{subject to}   & \Gamma x = Tz,& \\
                            & \l \leq x \leq u, & \\
                            & z \in \Z^{\mathscr{B}},
    \end{array}
\end{equation}
where $z \in \Z^\mathscr{B}$ is called cycle offset, an auxiliary integer variable akin to the periodic offset. 
In fact, given some feasible solution $(\pi, x, p)$ for the MIP \eqref{eq:pesp-formulation}, we can quickly recover the correct cycle offset $z$ for which the given $x$ is feasible by simply computing \begin{equation}\label{eq:p-to-z-offsets}
    z = \frac{\Gamma x}{T} = \frac{\Gamma(- B^\top \pi + T p)}{T} \overset{\text{\Cref{thm:cycle-periodicity}}}{=} \Gamma p.
\end{equation}
In this sense, we can also view \eqref{eq:pesp-formulation-cycle} as a variant of \eqref{eq:pesp-formulation}, where the periodic timetable variables $\pi$ have been eliminated and the linear independence between the periodic offset variables $p$ has been removed.
Although removed as variables, a coherently feasible timetable can easily be reconstructed from a given tension $x$ by fixing $\pi_v = 0$ for some vertex $v$, and then traversing the input graph in a depth-first fashion rooted at $v$ and assigning potentials according to the tension of the traversed arcs.

\begin{example}
    As the exemplary graph in \Cref{fig:ean} consists of a single cycle, its cycle space has rank 1 and the corresponding cycle matrix can be expressed as 
    \begin{equation}
        \Gamma =    \begin{pmatrix} 1 & -1 & 1 \end{pmatrix}.
    \end{equation}
    For the given periodic tension $x = (8,2,4)$ this results in $\Gamma x = 8-2+4 = 10$, i.e., $x$ corresponds to a cycle offset of $z = 1$ for period time $T=10$. 
\end{example}

The space of feasible periodic tensions is the same for the two MIP formulations \eqref{eq:pesp-formulation} and \eqref{eq:pesp-formulation-cycle}. 
In particular for the cycle formulation $X$ can be written as
\begin{equation}
    X \coloneqq \conv\left\{ x \in \R^{A(G)} \middle|\ \exists \ z \in \Z^{\mathscr{B}} \colon \Gamma x = Tz, \ \l \leq x \leq u \right\}
\end{equation}
However, we can as well consider the space of cycle offsets:
\begin{definition}\label{def:zonotope}
    Given a PESP instance $(G, T, \l, u, w)$ and an integral cycle basis $\mathscr{B}$ of $G$ with cycle matrix $\Gamma$, we define the \emph{cycle offset zonotope} as
        \begin{equation}
        Z \coloneqq \left\{ z \in \R^\mathscr{B} \ \middle|\ \exists x \in \R^{A(G)} \colon \Gamma x = T z, \ \l \leq x \leq u \right \}.
    \end{equation}
\end{definition}
We will see that $Z$ is indeed a zonotope in \Cref{sec:zonotopes}. Clearly, the integer points in $Z$ correspond to feasible cycle offset vectors for \eqref{eq:pesp-formulation-cycle}.

\section{The Tropical Tiling of the Periodic Timetable Space}
\label{sec:tropical}

Let $(G, T, \l, u, w)$ be a PESP instance with periodic timetable space $\Pi$. 

\begin{definition}\label{def:RofP}
    For $p \in \Z^{A(G)}$ we define
    \begin{equation}
        R(p) \coloneqq \left\{ \pi \in \R^{V(G)} \ \middle| \ \forall (i,j) \in A(G) \colon \l_{ij} -  T p_{ij} \leq \pi_j - \pi_i \leq u_{ij} -  T p_{ij}  \right\}.
    \end{equation}
\end{definition}
This definition leads to the decomposition $\Pi = \bigcup_{p \in \Z^{A(G)}} R(p)$ of the space of feasible periodic timetables into polyhedra $R(p)$. 
Since we assumed that $u - \l < T$, the polyhedra $R(p)$ are pairwise disjoint.
Each $R(p)$ has a very rigid structure that depends on the PESP input graph. 

\subsection{Weighted Digraph Polyhedra}

Focusing on the individual polyhedra $R(p)$, we will recall some of the properties of weighted digraph polyhedra along the lines of \cite{JoswigLoho}. 
These polyhedra have been called shortest path polyhedra in \cite{schrijverBook}. 
Unsurprisingly, such objects arise from weighted digraphs.
Here onwards we will always assume our graphs to have $n$ vertices, and we will make frequent and implicit use of some standard isomorphism $\R^{V(G)} \to \R^n$, fixed a priori.
\begin{definition}
    Given a digraph $D$ with weights $\kappa \colon A(D) \to \R$, the \emph{weighted digraph polyhedron} of $(D, \kappa)$, denoted by $W(D, \kappa)$, is defined by the points $\pi \in \R^n$ that satisfy the inequalities
    \begin{equation}
        \forall (i,j) \in A(D) \colon \pi_j - \pi_i \leq \kappa_{ij}.
    \end{equation}
\end{definition}
The points of a weighted digraph polyhedron correspond to \emph{feasible potentials} of the weighted digraph, again see \cite{schrijverBook}.
In particular, $W(D, \kappa)$ is empty if and only if $\kappa$ is not conservative, i.e., there exists a negative cycle in $(D, \kappa)$.

A general observation is that any feasible potential of $(D, \kappa)$, i.e., any point of $W(D, \kappa)$, can be translated by any real multiple of the all-ones vector $\mathbf{1} = (1,\ldots,1)$ and it remains feasible. 
Indeed, if $\pi_j - \pi_i \leq \kappa_{ij}$ then $(\pi_j + k) - (\pi_i + k) \leq \kappa_{ij}$ for all $k \in \R$. 
This implies that any non-empty weighted digraph polyhedron contains $\R\mathbf{1} = \langle \mathbf{1} \rangle$ in its lineality space.
In fact, if $D$ is weakly connected and $W(D, \kappa)$ is non-empty, then the lineality space of $W(D, \kappa)$ equals $\R\mathbf{1}$.
Moreover, if $D$ is also strongly connected, then the recession cone of $W(D, \kappa)$ is only its lineality space \cite[\S2.3-\S2.4]{JoswigLoho}.

\begin{figure}[htbp]
    \centering
    	\begin{tikzpicture}[scale=2.5]
		\tikzstyle{p} = [line width=1.5, ->]
		\tikzstyle{u} = [line width=1.5]
		\tikzstyle{q} = [line width=1.5, ->]
		\tikzstyle{r} = [line width=1.5, ->]
		\tikzstyle{v} = [draw, circle, inner sep=1, minimum width=12, font=\footnotesize]
		\tikzstyle{w} = [v]
		\tikzstyle{s} = [draw, rectangle, inner sep=1, minimum width=12, minimum height=12, font=\footnotesize]
		\tikzstyle{t} = [midway, font=\footnotesize, sloped]
		\tikzstyle{ta} = [t, above]
		\tikzstyle{tb} = [t, below]
		\tikzstyle{tr} = [t, right]
		\tikzstyle{tl} = [t, left]
		\tikzstyle{l} = [font=\scriptsize]
		\node[w] (A) at (0, 0) {$v_0$ };
		\node[v] (B) at (1, 1.5) {$v_1$};
		\node[w] (C) at (2, 0) { $v_2$};
		
		\draw[p] (A) edge[looseness=1,out=110,in=180] node[ta] {$12$} (B);
		\draw[p] (B) edge[looseness=1,out=0,in=70] node[ta] {$3$} (C);
		\draw[p] (A) edge[looseness=1,out=-45,in=-135] node[tb] {$10$} (C);	
		
		\draw[p] (B) -- node[ta] {$-3$} (A);
		\draw[p] (C) -- node[ta,black] {$6$} (B);
		\draw[p] (C) -- node[ta,black] {$-2$} (A);

	\end{tikzpicture}
    \caption{
        Weighted digraph $(\overline{G}, \kappa(p))$ for $p = (0,0,1)$ of $G$ in \Cref{fig:ean}.
        Edge labels are $\kappa(p)_{ij}$.
        }
    \label{fig:goverline}
\end{figure}
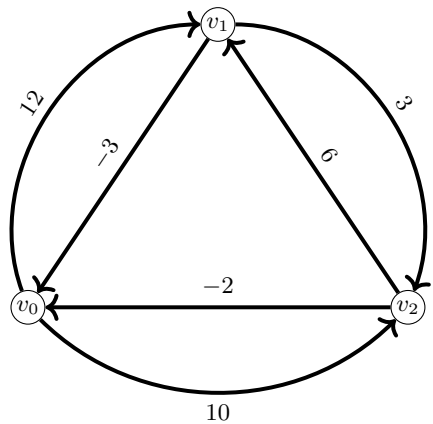

In our particular case of periodic timetabling, we construct for each $p \in \Z^{A(G)}$ a weighted digraph $(\overline{G}, \kappa(p))$ whose weighted digraph polyhedron is $R(p)$. 
The vertex set of $\overline{G}$ is $V(\overline{G}) \coloneqq V(G)$, and its arc set is $A(\overline{G}) \coloneqq A(G) \cup A(G^\top)$, where $A(G^\top) = \{ (j, i) | (i,j) \in A(G)\}$.
Finally the arc weights are $\kappa(p)_{ij} \coloneqq u_{ij} - T p_{ij}$ for all $(i,j) \in A(G)$ and $\kappa(p)_{ij} \coloneqq T p_{ji} - \l_{ji}$ for all $(i,j) \in A(G^\top)$.
By construction every arc of $\overline{G}$ will always have a corresponding antiparallel arc, so we understand that the resulting inequalities consist in pairs of parallel halfspaces with opposing directions.
An example of the $\overline{G}$ construction is seen in \Cref{fig:goverline}.

\subsection{Polytropes}

A different point of view is that of tropical geometry, which is algebraic geometry over the \emph{tropical semiring} $\T = \left(\R \cup \{\infty\}, \oplus, \odot\right)$.
For any two $a,b \in \T$ the tropical operations are defined as 
\begin{equation}
    a \oplus b \coloneqq \min\{a,b\} \quad \text{~and~} \quad a \odot b \coloneqq a + b .
\end{equation}
The tropical semiring has $\infty$ and $0$ as additive and multiplicative identity elements, respectively.
Lastly, as no element in $\T$ except $\infty$ has an additive inverse, we can understand tropical numbers to have no sign in the sense we may be used to.
What mainly interests us is the theory of so-called tropical convexity, whose definition is similar to the traditional setting.
\begin{definition}
    A set $S \subset \T^n$ is \emph{tropically convex} if 
    \begin{equation}
        (a\odot x) \oplus (b\odot y) \in S
    \end{equation}
    for every $x,y \in S$, and any $a,b \in \T$.    
\end{definition}
No restrictions are posed on the tropical scalar coefficients, which implies that any tropically convex set is closed under tropical scalar multiplication.
Alternatively this means that non-empty tropically convex sets are unbounded in $\T^n$ as they must always contain $\T \odot (0, \ldots, 0) = \R \mathbf{1} \cup \{(\infty, \ldots, \infty)\}$.
We will therefore consider the \emph{tropical projective space} $\TP^{n-1}$, which is defined as the quotient of $(\T^n \setminus \{(\infty, \ldots, \infty)\}) / \T$ via tropical scalar multiplication. 
The tropical projective space $\TP^{n-1}$ contains the affine chart $\R^n/\R\mathbf{1} \cong \R^{n-1}$, which is called \emph{tropical projective $(n-1)$-torus} in \cite{ETC}, or \emph{tropical affine space} in \cite{JoswigKulas}. 
Other authors sometimes have referred to $\R^n/\R\mathbf{1} \cong \R^{n-1}$ also as tropical projective space \cite{Develin2004}. 
We point to \cite{Develin2004} and \cite{hill2020tropical} for a thorough introduction to tropical convexity, and \cite{ETC} for a general compendium of tropical geometry and combinatorics.

It is interesting to point out that tropical convexity and traditional convexity do not imply each other, which draws attention to the cases where, instead, both properties are upheld.
\begin{definition}
    A \emph{polytrope} is the tropical convex hull of a finite set of points in $\R^n/\R\mathbf{1}$ that is also traditionally convex.
\end{definition}
There is a correspondence between polytropes and weighted digraph polyhedra:
\begin{theorem}[\cite{JoswigLoho}, Proposition 48]
    If $(D, \kappa)$ is a strongly connected weighted digraph on $n$ vertices, then the quotient modulo $\R\mathbf{1}$ of its weighted digraph polyhedron $W(D, \kappa) \subseteq \R^n$ is a polytrope in $\R^n/\R\mathbf{1}$. 
    Moreover, every polytrope in $\R^n/\R\mathbf{1}$ arises this way.
\end{theorem}
The minimal generating set of a polytrope in $\R^n/\R\mathbf{1}$ is a finite set of points, which are called the \emph{tropical vertices} of the polytrope.
Such set is known to have fixed size $n$.
This interestingly means that all polytropes are tropical simplices in tropical projective space.

Returning to the space of feasible periodic timetables of a PESP instance, we obtain with \Cref{def:RofP} and the previous result 
\begin{theorem}
    The space $\Pi/\R\mathbf{1} \subset \TP^{n-1}$ decomposes into the union of the pairwise disjoint polytropes $R(p)/\R\mathbf{1}$.
\end{theorem}

\subsection{The Torus of Feasible Periodic Timetables}

In this section we will always be working with the PESP instance $(G, T, \l, u, w)$, where $|V(G)| = n$.
One can easily observe that there are infinitely many disjoint polytropes $R(p)/\R\mathbf{1}$ which partition $\Pi/\R\mathbf{1}$.
However, many of them are simply translated copies of each other. 
This is due to the following symmetry condition, which is characteristic of PESP itself: 
\begin{lemma}\label{lemma:piPlusTq-symmetry}
    If $\pi \in \Pi$ then $\pi' \coloneqq \pi + Tq \in \Pi$ for all $q \in \Z^{V(G)}$.
\end{lemma}
\begin{proof}
    Let $p \in \Z^{A(G)}$ be the periodic offset such that $\l_{ij} \leq \pi_j - \pi_i + Tp_{ij} \leq u_{ij}$. 
    Then we can construct $p'$ such that $p'_{ij} \coloneqq  p_{ij} - q_j + q_i$.
    Then $\l_{ij} \leq \pi'_j - \pi'_i + Tp'_{ij} \leq u_{ij}$ for all $(i,j) \in A(G)$.
\end{proof}
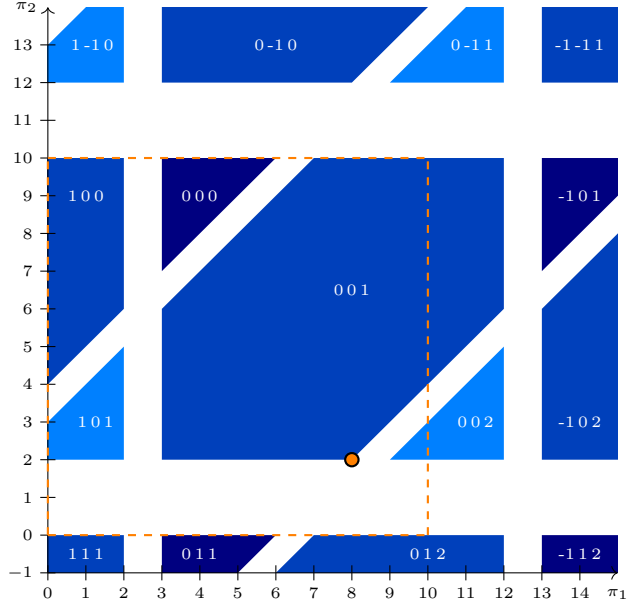
\begin{figure}
    \centering
    	\begin{tikzpicture}[scale=.5]
	\definecolor{color0}{HTML}{000080}
	\definecolor{color1}{HTML}{0040BB}
	\definecolor{color2}{HTML}{0080FF}
	\tikzstyle{l} = [draw, line width=0.8, grau]
	\tikzstyle{v} = [circle, draw, black, line width=0.8, fill=white, inner sep=1.5]
	\tikzstyle{f} = []
	\tikzstyle{t} = [white, font=\tiny,  ]
	\tikzstyle{g} = [orange, dashed, line width=0.8]
	\draw[->] (0, -1) -- (15, -1);
	\draw[->] (0, -1) -- (0, 14);

	\foreach \i in {0,1,...,7,8, 9,10,...,14} {
		\draw (\i,-1.2) -- (\i,-0.8) node[below,yshift=-1ex,font=\tiny] {$\i$};
	};
	\foreach \i in {-1,0,1,2,3,4,...,13} {
		\draw (-0.2,\i) -- (0.2,\i) node[left,xshift=-1ex,font=\tiny] {$\i$};
	};
	
	\node[anchor=north, yshift=-0.5ex, font=\tiny] at (15,-1) {$\pi_1$};
	\node[anchor=east, font=\tiny] at (0,14) {$\pi_2$};
	\begin{pgfonlayer}{bg}
	\clip (0,-1) rectangle (15,14);
	\fill[f, fill=color0] (3,7) -- (3,10) -- (6,10) -- cycle;
	\node[t, anchor=west] at (3.25, 9) {0\,0\,0};
	\fill[f, fill=color0] (13,7) -- (13,10) -- (16,10) -- cycle;
	\node[t] at (14, 9) {-1\,0\,1};
	\fill[f, fill=color0] (3,-3) -- (3,0) -- (6,0) -- cycle;
	\node[t, anchor=west] at (3.25, -.5) {0\,1\,1};
	\fill[f, fill=color0] (13,-3) -- (13,0) -- (16,0) -- cycle;
	\node[t] at (14, -.5) {-1\,1\,2};
	
	\fill[f, fill=color1] (3,2) -- (3,6) -- (7,10) -- (12,10) -- (12,6) -- (8,2) -- cycle;
	\node[t] at (8, 6.5) {0\,0\,1};
	\fill[f, fill=color1] (13,2) -- (13,6) -- (17,10) -- (22,10) -- (22,6) -- (18,2) -- cycle;
	\node[t] at (14, 3) {-1\,0\,2};
	\fill[f, fill=color1] (-7,2) -- (-7,6) -- (-3,10) -- (2,10) -- (2,6) -- (-2,2) -- cycle;
	\node[t] at (1, 9) {1\,0\,0};
	
	\fill[f, fill=color1] (3,12) -- (3,16) -- (7,20) -- (12,20) -- (12,16) -- (8,12) -- cycle;
	\node[t] at (6, 13) {0\,-1\,0};
	\fill[f, fill=color1] (13,12) -- (13,16) -- (17,20) -- (22,20) -- (22,16) -- (18,12) -- cycle;
	\node[t] at (14, 13) {-1\,-1\,1};
	
	\fill[f, fill=color1] (3,-8) -- (3,-4) -- (7,0) -- (12,0) -- (12,-4) -- (8,-8) -- cycle;
	\node[t] at (10, -0.5) {0\,1\,2};
	\fill[f, fill=color1] (-7,-8) -- (-7,-4) -- (-3,0) -- (2,0) -- (2,-4) -- (-2,-8) -- cycle;
	\node[t] at (1, -0.5) {1\,1\,1};
	
	\fill[f, fill=color2] (9,2) -- (12,5) -- (12,2) -- cycle;
	\node[t, anchor=east] at (12, 3) {0\,0\,2};
	\fill[f, fill=color2] (-1,2) -- (2,5) -- (2,2) -- cycle;
	\node[t, anchor=east] at (2, 3) {1\,0\,1};
	
	\fill[f, fill=color2] (9,12) -- (12,15) -- (12,12) -- cycle;
	\node[t, anchor=east] at (12, 13) {0\,-1\,1};
	\fill[f, fill=color2] (-1,12) -- (2,15) -- (2,12) -- cycle;
	\node[t, anchor=east] at (2, 13) {1\,-1\,0};

	\fill[circle, draw, line width=0.8, fill=orange, inner sep=2] (8,2) circle (5 pt);
	\end{pgfonlayer}

	\draw[g] (0,0) -- (10,0);%
	\draw[g] (0,10) -- (10,10);
	\draw[g] (0,0) -- (0,10);
	\draw[g] (10,0) -- (10,10);
\end{tikzpicture}
    \caption{
        Space of potentials $\Pi/\R\mathbf{1}$ and its decomposition into polytropes of \Cref{fig:ean}. 
        Each region is labeled by the periodic offset vector $p = (p_{01}, p_{02}, p_{12})$. 
        The marked corner corresponds to the labeled timetable of \Cref{fig:ean}.
        }
    \label{fig:tiling}
\end{figure}
As a consequence we can reduce our feasible space even further, by taking another quotient.
The intention is to incorporate the equivalence relation $\pi' \cong \pi + T q$ for all $q \in \Z^{V(G)}$, and the resulting space is therefore 
\begin{equation}\label{eq:torusSpace}
    \mathscr{T} \coloneqq \left(\R^{V(G)}/\left(T\Z\right)^{V(G)}\right)/\R\mathbf{1},
\end{equation}
namely a $(n-1)$-dimensional torus of side length $T$.
It is therefore sufficient in $\R^n/\R\mathbf{1}$ to consider any $(n-1)$-hypercube of side length $T$, positioned anywhere.
Any such hypercube would in fact be a representative of $\mathscr{T}$, and we will call said region a fundamental domain.
Going forward we will also abbreviate the following $(R(p)/(T\Z)^n)/\R\mathbf{1} \eqqcolon \mathbf{R}(p) \subseteq \mathscr{T}$.

\begin{lemma}
\label{lem:RpEquiv}
    Let $p, p' \in \Z^A$. 
    Then $\mathbf{R}(p) = \mathbf{R}(p')$ if and only if $\Gamma p = \Gamma p'$. 
\end{lemma}
\begin{proof}
    Let $\Gamma$ be the matrix of an integral cycle basis, for example a fundamental cycle basis, and $B$ the incidence matrix of $G$. 
    Suppose $\mathbf{R}(p) = \mathbf{R}(p') \subseteq \mathscr{T}$. 
    By \Cref{lemma:piPlusTq-symmetry}, we have $p' = p - B^\top q$. 
    Since $\Gamma B^\top = 0$ we obtain $\Gamma p' = \Gamma (p - B^\top q) = \Gamma p$.
    Conversely, if $\Gamma p = \Gamma p'$ then $p - p' \in \ker\Gamma$. 
    By \Cref{thm:cycle-periodicity} then $\im B^\top = \ker \Gamma$, and as the basis is chosen to be integral, there exists a $q \in \Z^V$ such that $p - p' = B^\top q$, and so $p' = p - B^\top q$. 
    Again by \Cref{lemma:piPlusTq-symmetry}, we have $\mathbf{R}(p) = \mathbf{R}(p')$.
\end{proof}
It is clear from this proof that any two periodic offset vectors $p$ and $p'$ yield the same cycle offset $z$, and thereby the same polytrope $\mathbf{R}(p) \subseteq \mathscr{T}$, whenever they differ by some linear combination of the rows of $B$, since $\im B^\top = \ker \Gamma$.
This motivates us to define $\mathfrak{R}(z) \coloneqq \mathbf{R}(p)$ for any $p$ such that $\Gamma p = z$.

\begin{figure}
    \centering
    \includegraphics[width=0.6\textwidth]{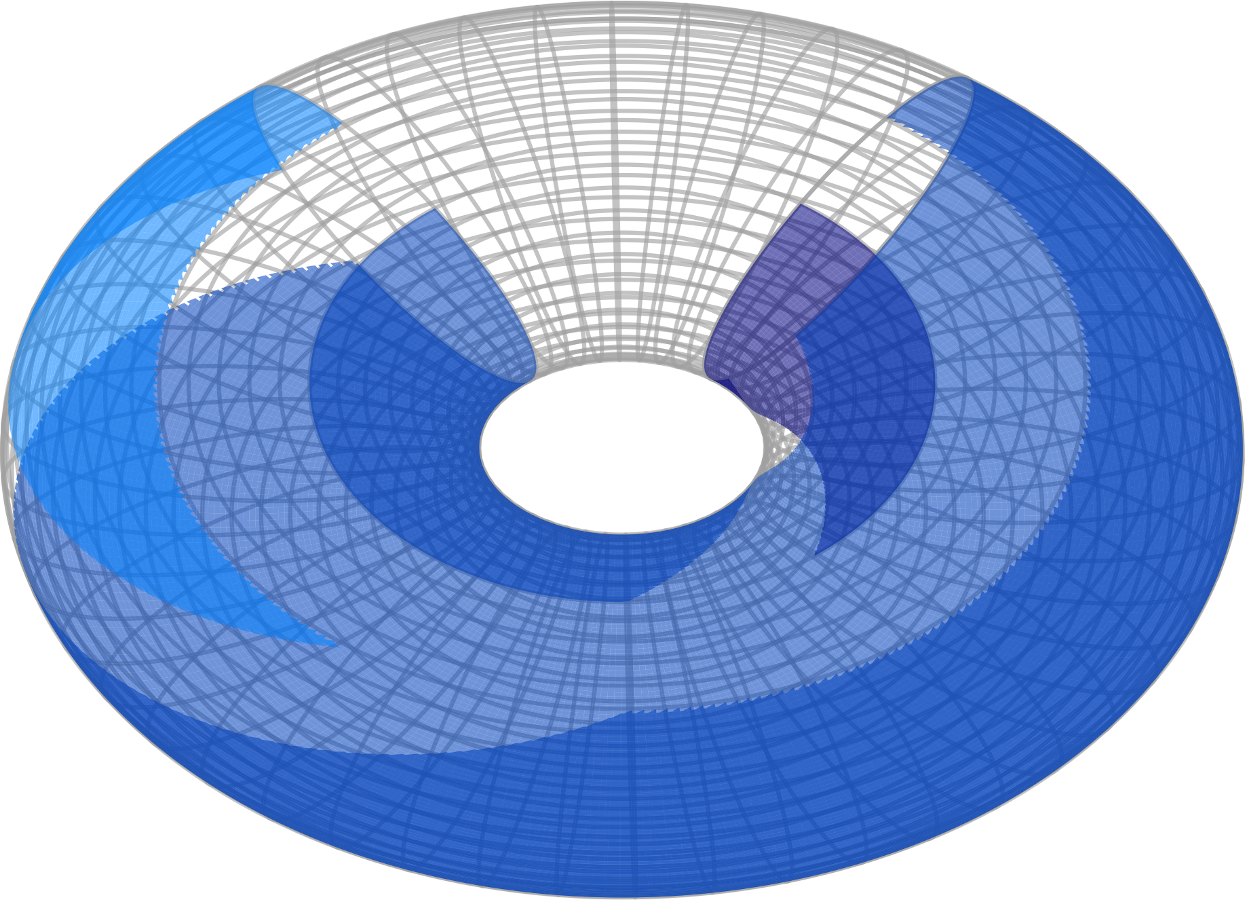}
    \caption{
        The orange fundamental domain shown in \Cref{fig:tiling} is wrapped onto the torus $\mathscr{T}$.
        }
    \label{fig:torus}
\end{figure}

\begin{example}\label{ex:polytropalTiling} 
    Again, we can visualise properties of \Cref{lemma:piPlusTq-symmetry} and \Cref{lem:RpEquiv} with the help of our running example from \Cref{fig:tiling}.
    Clearly, we have $T$-translates of the same regions in $\Pi/\R\mathbf{1}$ as indicated by the different shades of blue.
    Using cycle matrix $\Gamma = (1, -1, 1)$, we can, indeed, identify the light blue triangular polytropes in $\Pi/\mathbb{R}\mathbf{1}$ with the cycle offset $z = \Gamma p = 2$, the dark blue ones with $z = 0$, and the hexagonal polytropes with $z = 1$ (also recall the graph of \Cref{fig:goverline}, which generates the depicted hexagon as a weighted digraph polyhedron). 
    This allows us to restrict ourselves to only three regions in the torus $\mathscr{T}$.
    A possible choice for a fundamental domain is visualised by the orange dashed square.
    The space depicted in \Cref{fig:tiling} is then quotiented as described in \Cref{eq:torusSpace} and the result is depicted in \Cref{fig:torus}.
\end{example}

\begin{lemma}\label{lemma:SymmetricPiSameX}
    Let $\pi, \pi' \in \Pi/\R\mathbf{1}$ be two periodic timetables, and let $x \in \R^{A(G)}$ be the feasible periodic tension associated to $\pi$.
    Then there exists $q \in \Z^{V(G)}$ such that $\pi' = \pi + Tq$ if and only if $x$ is the feasible periodic tension associated to $\pi'$ as well.
\end{lemma}
\begin{proof}
    Supposing there is such a $q$, then by definition we have 
    \begin{equation} 
        \begin{split}
            x_{ij}  &= \pi_j - \pi_i + Tp_{ij} = \pi_j + Tq_j - \pi_i - Tq_i + Tp_{ij} - Tq_j + Tq_i \\ 
                    &= \pi'_j - \pi'_i + Tp'_{ij}
        \end{split}
    \end{equation}
    for every arc $(i,j) \in A(G)$.
    This way we implicitly constructed the appropriate periodic offset $p'$, through which $x$ is indeed the feasible periodic tension associated to $\pi'$.
    
    Conversely, suppose $x$ is such a tension for $\pi'$ as well as $\pi$, i.e.,  
    \begin{equation}
        -B^\top \pi + Tp = -B^\top \pi' + Tp' 
    \end{equation}
    holds for appropriate $p, p'$. 
    An application of the cycle matrix from the left yields 
    \begin{equation}
        \Gamma p = \Gamma p',
    \end{equation}
    since $\Gamma B^\top = 0$.
    The claim then follows directly from \Cref{lemma:piPlusTq-symmetry} and \Cref{lem:RpEquiv}.
\end{proof}

\subsection{Dimension of Timetable Polytropes}

Focusing on the dimension of a specific $\mathbf{R}(p)$, we have already dealt with the case of dimension $-1$, i.e., when $(\overline{G}, \kappa(p))$ has a cycle of negative weight.
Otherwise we can compute $\dim \mathbf{R}(p)$ following \cite{JoswigLoho}.
The equality graph of $(\overline{G}, \kappa(p))$, denoted by $E(\overline{G}, \kappa(p))$, is an undirected graph on the same vertex set $V(G)$ such that two nodes are adjacent if and only if they are contained in a cycle of weight 0 in $(\overline{G}, \kappa(p))$.
\begin{theorem}[\cite{JoswigLoho}, Lemma~5]
    If $\mathbf{R}(p) \neq \varnothing$, then $\dim \mathbf{R}(p)$ is equal to the number of connected components of $E(\overline{G}, \kappa(p))$ minus 1.
\end{theorem}
The more the cycles of weight 0 are either numerous or large, the lower the dimension of $\mathbf{R}(p)$.
It is a single point in the torus if and only if all vertices in $V(G)$ are contained in the same cycle of weight 0.
Geometrically this phenomenon is easily understood.
Supposing we had a cycle of weight 0 with $k$ arcs, then any $k-1$ of them define hyperplanes in the construction of $\mathbf{R}(p)$ as a weighted digraph polyhedron, and their intersection is a $n-k-1$ linear variety $L \subset \mathscr{T}$.
The intersection of the corresponding halfspaces is then, in any case, some full dimensional orthant non-differentiable in $L$.
The last arc of the cycle now defines a hyperplane that, again, goes through $L$, and the direction of the corresponding halfspace is outside of the orthant, thereby restricting $\mathbf{R}(p) \subseteq L$.
In this light we have a new way of understanding how negative cycles immediately make polytropes empty, since cycles of weight 0 are in fact limit cases of such behaviour.
In general we can also a priori predict the existence of not full-dimensional $\mathbf{R}(p) \subseteq \mathscr{T}$ for certain vectors $p$, by inspecting the weighted digraph $(\overline{G}, \kappa(\cdot))$.
For every cycle $\gamma$ of weight $kT + \gamma \cdot p$, with $k \in \Z$, there will exist appropriate periodic offsets $p$ so as for $\gamma$ to have weight 0 in $(\overline{G}, \kappa(p))$.

\subsection{Embedding Timetable Polytropes into the Periodic Tension Polytope}

Up until now, we mainly focused on the timetable space $\Pi$.
As the two equivalent formulations of PESP suggest, one can also look at the same problem from a different perspective, namely by considering the tension space.
Recall that periodic tension polytope is defined, using cycle offset variables, as
\begin{equation}
    X = \conv\left\{ x \in \R^{A(G)} \middle|\ \exists \ z \in \Z^{\mathscr{B}} \colon \Gamma x = Tz, \ \l \leq x \leq u \right\},
\end{equation}
for some cycle basis $\mathscr{B}$ of $G$, with corresponding matrix $\Gamma$.
There is a natural map from $R(p)/\R \mathbf{1}$ to $X$. 
More precisely, we can state the following: 
\begin{lemma}
    For a non-empty $R(p)$ of the PESP instance $(G, T, \l, u, w)$, the affine map \begin{equation}
        \begin{split}
            m_p :   R(p)/\R \mathbf{1} &\longrightarrow X \\ 
                    \pi &\longmapsto -B^\top\pi + Tp
        \end{split}
    \end{equation}
    embeds $R(p)/\R \mathbf{1}$ in $X$. 
\end{lemma}
This clearly extends to $\mathscr{T}$ as well:
\begin{corollary}
    For a region $\mathfrak{R}(z) \in \mathscr{T}$ the affine map
    \begin{equation}
        \begin{split}
        \mathfrak{m}_z :   \mathfrak{R}(z) &\longrightarrow X \\ 
                \pi &\longmapsto -B^\top\pi + Tp
        \end{split}
    \end{equation}
    for some $p$ such that $z = \Gamma p$ embeds $\mathfrak{R}(z)$
    in $X$.
\end{corollary}

\begin{figure}[htbp]
    \centering
    \begin{tikzpicture}[
x  = {(0.9cm,0.25cm)},
y  = {(0cm,1cm)},
z  = {(0.75cm,-0.4cm)},
    x = {(0.7cm,0.4cm)},
    y = {(0cm,1cm)},
    z = {(0.9cm,-0.2cm)},
scale = 0.4,
color = {lightgray}]

\definecolor{pointcolor_r}{rgb}{ 1,0,0 }
\tikzstyle{pointstyle_r} = [fill=pointcolor_r]


\coordinate (v2_r) at (3, 7, 4); 
\coordinate (v3_r) at (3, 10, 7);
\coordinate (v4_r) at (6, 10, 4);

\coordinate (v7_r) at (9, 2, 13);
\coordinate (v8_r) at (12, 2, 10);
\coordinate (v9_r) at (12, 5, 13);

\coordinate (v10_r) at (12, 6, 4);
\coordinate (v11_r) at (12, 10, 8);
\coordinate (v5_r) at (8, 2, 4); 
\coordinate (v6_r) at (7, 10, 13); 
\coordinate (v0_r) at (3, 2, 9); 
\coordinate (v1_r) at (3, 6, 13);

\definecolor{edgecolor_r}{rgb}{ 0,0,0 }

\definecolor{gcolor}{rgb}{ 0.3667,0.9255,0.5196 }
\definecolor{rcolor}{rgb}{ 0.9255,0.5196,0.3667 }
\definecolor{bcolor}{rgb}{ 0.5196,0.3667,0.9255 }
\definecolor{scolor}{rgb}{ 0.925,0.925,0.925 }
\definecolor{color0}{HTML}{000080}
\definecolor{color1}{HTML}{0040BB}
\definecolor{color2}{HTML}{0080FF}

\tikzstyle{facestyle_r} = [fill=scolor, fill opacity=0.4, draw=edgecolor_r, line width=1 pt, line cap=round, line join=round]
\tikzstyle{bound} = [facestyle_r]
\tikzstyle{cycle} = [fill opacity=1, draw=edgecolor_r]
\tikzstyle{ccycle} = [facestyle_r]
\tikzstyle{vertex} = [text=black, inner sep=0.5pt, align=left, draw=none, scale = 0.8] 

%



\draw[thick,->] (0,0,0) -- (3,0,0) node[anchor=south]{$x_{01}$};
    \draw[thick,->] (0,0,0) -- (0,3,0) node[anchor=north west]{$x_{02}$};
    \draw[thick,->] (0,0,0) -- (0,0,3) node[anchor=north]{$x_{12}$};

\draw[facestyle_r, color1, cycle] (v6_r) -- (v1_r) -- (v0_r) -- (v5_r) -- (v10_r) -- (v11_r) -- cycle;

\draw[facestyle_r, bound] (v10_r) -- (v5_r) -- (v2_r) -- (v4_r) -- (v10_r) -- cycle;
\draw[facestyle_r] (v8_r) -- (v5_r) -- (v10_r) -- (v8_r) -- cycle;
\draw[facestyle_r, bound] (v0_r) -- (v5_r) -- (v8_r) -- (v7_r) -- (v0_r) -- cycle;

\draw[facestyle_r, color0, cycle] (v4_r) -- (v2_r) -- (v3_r) -- (v4_r) -- cycle;
\draw[facestyle_r, ccycle] (v2_r) -- (v5_r) -- (v0_r) -- (v2_r) -- cycle; 
\draw[facestyle_r] (v0_r) -- (v7_r) -- (v1_r) -- (v0_r) -- cycle; 
\draw[facestyle_r, color2, cycle] (v7_r) -- (v8_r) -- (v9_r) -- (v7_r) -- cycle;

\draw[facestyle_r] (v11_r) -- (v6_r) -- (v9_r) -- (v11_r) -- cycle;
\draw[facestyle_r] (v10_r) -- (v4_r) -- (v11_r) -- (v10_r) -- cycle;
\draw[facestyle_r, bound] (v10_r) -- (v11_r) -- (v9_r) -- (v8_r) -- (v10_r) -- cycle;
\draw[facestyle_r, bound] (v7_r) -- (v9_r) -- (v6_r) -- (v1_r) -- (v7_r) -- cycle;
\draw[facestyle_r, bound] (v2_r) -- (v0_r) -- (v1_r) -- (v3_r) -- (v2_r) -- cycle;
\draw[facestyle_r, bound] (v4_r) -- (v3_r) -- (v6_r) -- (v11_r) -- (v4_r) -- cycle;
\draw[facestyle_r] (v3_r) -- (v1_r) -- (v6_r) -- (v3_r) -- cycle;

\fill[pointcolor_r] (v2_r) circle (1 pt);
\node[vertex, right] at (v2_r) {$(3,7,4)$};
\fill[circle, draw, black, line width=0.8, fill=orange, inner sep=2] (v5_r) circle (5 pt);
\node[vertex, above right] at (v5_r)  {$(8,2,4)$};
\fill[pointcolor_r] (v4_r) circle (1 pt);
\node[vertex, above right] at (v4_r) {$(6,10,4)$};
\fill[pointcolor_r] (v0_r) circle (1 pt);
\node[vertex, below right] at (v0_r)  {$(3,2,9)$};
\fill[pointcolor_r] (v10_r) circle (1 pt);
\node[vertex, above right] at (v10_r)  {$(12,6,4)$};
\fill[pointcolor_r] (v3_r) circle (1 pt);
\node[vertex, above right] at (v3_r){$(3,10,7)$};
\fill[pointcolor_r] (v1_r) circle (1 pt);
\node[vertex, below right] at (v1_r) {$(3,6,13)$};
\fill[pointcolor_r] (v7_r) circle (1 pt);
\node[vertex, below right] at (v7_r) {$(9,2,13)$};
\fill[pointcolor_r] (v8_r) circle (1 pt);
\node[vertex, above left] at (v8_r) {$(12,2,10)$};
\fill[pointcolor_r] (v11_r) circle (1 pt);
\node[vertex, above right] at (v11_r) {$(12,10,8)$};
\fill[pointcolor_r] (v6_r) circle (1 pt);
\node[vertex, above right] at (v6_r){$(7,10,13)$};
\fill[pointcolor_r] (v9_r) circle (1 pt);
\node[vertex, above right] at (v9_r) {$(12,5,13)$};


\end{tikzpicture}
    \caption{Periodic tension polytope $X$ with the embeddings of $\mathfrak{R}(z)$ for $z \in \{0,1,2\}$, corresponding to the instance of \Cref{ex:running} and its tiling, shown in \Cref{fig:tiling}.}
    \label{fig:tensionPolytope}
\end{figure}
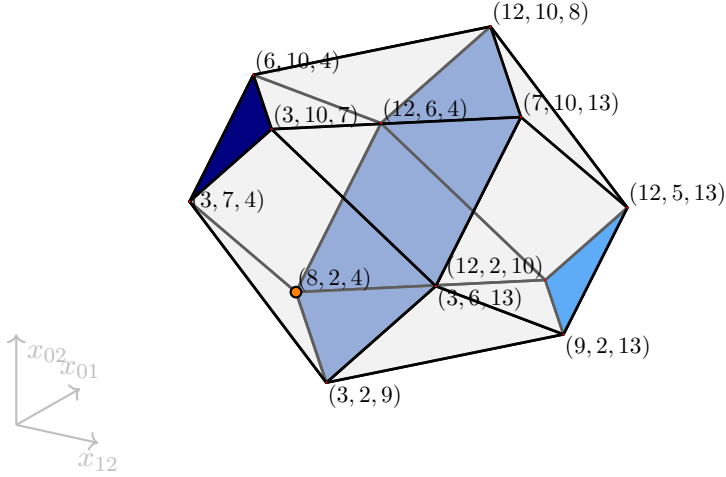

A visualization of these two properties can be seen in \Cref{fig:tensionPolytope}. Here, $\mathfrak{m}_1$ maps the hexagonal polytrope $\mathfrak{R}(1) = \mathbf{R}((0,0,1))$ of the tiling in \Cref{fig:tiling} to the hexagon in the interior of $X$, while the two triangular polytropes corresponding to $\mathfrak{R}(0) = \mathbf{R}((0,0,0))$ and $\mathfrak{R}(2) = \mathbf{R}((0,0,2))$, are mapped to facets of the tension polytope. 

While these observations might seem trivial at first glance, they allow for a different geometric interpretation. 
Consider the \emph{periodic tension polytope of the LP-relaxation} of a PESP instance, i.e., 
\begin{equation}
    X_\text{LP} \coloneqq \conv\left\{ x \in \R^{A(G)} \middle|\ \exists \ z \in \R^{\mathscr{B}} \colon \Gamma x = Tz, \ \l \leq x \leq u \right\}.
\end{equation}
We may also call this object the \emph{fractional periodic tension polytope}.
As an appropriately real-valued $z$ always exists, it is easy to see that $X_\text{LP}$ simply corresponds to a hyperrectangle with side lengths given by the span of the arcs.
More precisely 
\begin{equation}
    X_\text{LP} = \bigtimes_{(i,j)\in A(G)} [l_{ij}, u_{ij}].
\end{equation}
Going forward, we will for simplicity say cube rather than hyperrectangle, especially because most of our interest will be in the combinatorial properties of the faces of the object.

\begin{figure}[htbp]
    \centering
    \begin{tikzpicture}[
x  = {(0.9cm,0.25cm)},
y  = {(0cm,1cm)},
z  = {(0.75cm,-0.4cm)},
    x = {(0.7cm,0.4cm)},
    y = {(0cm,1cm)},
    z = {(0.9cm,-0.2cm)},
scale = 0.4,
color = {lightgray}]

\definecolor{pointcolor_r}{rgb}{ 1,0,0 }
\tikzstyle{pointstyle_r} = [fill=pointcolor_r]


\coordinate (v2_r) at (3, 7, 4); 
\coordinate (v3_r) at (3, 10, 7);
\coordinate (v4_r) at (6, 10, 4);

\coordinate (v7_r) at (9, 2, 13);
\coordinate (v8_r) at (12, 2, 10);
\coordinate (v9_r) at (12, 5, 13);

\coordinate (v10_r) at (12, 6, 4);
\coordinate (v11_r) at (12, 10, 8);
\coordinate (v5_r) at (8, 2, 4); 
\coordinate (v6_r) at (7, 10, 13); 
\coordinate (v0_r) at (3, 2, 9); 
\coordinate (v1_r) at (3, 6, 13); 

\coordinate (c1) at (3,2,4);
\coordinate (c2) at (3,2,13);
\coordinate (c3) at (3,10, 13);
\coordinate (c4) at (3, 10, 4); 
\coordinate (c5) at (12, 2, 4);
\coordinate (c6) at (12,2,13);
\coordinate (c7) at (12,10, 13);
\coordinate (c8) at (12, 10, 4); 

\coordinate (h21) at (6, 14, 8);
\coordinate (h22) at (11, 10, -1); 
\coordinate (h24) at (3, 3, 0);
\coordinate (h23) at (-2, 7, 9); 


\coordinate (h11) at (11,14,13);
\coordinate (h12) at (16, 10, 4); 
\coordinate (h14) at (4, -2, 4);
\coordinate (h13) at (-1, 2, 13); 

\coordinate (h31) at (12,9,17);
\coordinate (h32) at (17, 5, 8); 
\coordinate (h34) at (9, -2, 9);
\coordinate (h33) at (4, 2, 18); 

\definecolor{edgecolor_r}{rgb}{ 0,0,0 }
\definecolor{gcolor}{rgb}{ 0.3667,0.9255,0.5196 }
\definecolor{rcolor}{rgb}{ 0.9255,0.5196,0.3667 }
\definecolor{bcolor}{rgb}{ 0.5196,0.3667,0.9255 }
\definecolor{scolor}{rgb}{ 0.925,0.925,0.925 }
\definecolor{color0}{HTML}{000080}
\definecolor{color1}{HTML}{0040BB}
\definecolor{color2}{HTML}{0080FF}

\tikzstyle{facestyle_r} = [fill=scolor, fill opacity=0.4, draw=edgecolor_r, line width=1 pt, line cap=round, line join=round]
\tikzstyle{bound} = [facestyle_r]
\tikzstyle{cycle} = [fill opacity=1, draw=edgecolor_r]
\tikzstyle{ccycle} = [facestyle_r]
\tikzstyle{vertex} = [text=black, inner sep=0.5pt, align=left, draw=none, scale = 0.8] 

%



\draw[facestyle_r] (c5) -- (c6) -- (c7) -- (c8) -- cycle;
\draw[facestyle_r] (c1) -- (c2) -- (c6) -- (c5) -- cycle;
\draw[facestyle_r] (c1) -- (c4) -- (c8) -- (c5) -- cycle;

\draw[facestyle_r, color0!30]  (h21) -- (h22) -- (h24) -- (h23) -- cycle; 
\draw[facestyle_r, color1!30]  (h11) -- (h12) -- (h14) -- (h13) -- cycle; 
\draw[facestyle_r, color2!30]  (h31) -- (h32) -- (h34) -- (h33) -- cycle; 

\draw[facestyle_r, color1, opacity = .5] (v6_r) -- (v1_r) -- (v0_r) -- (v5_r) -- (v10_r) -- (v11_r) -- cycle;
\draw[facestyle_r, color0, opacity = .5] (v4_r) -- (v2_r) -- (v3_r) -- (v4_r) -- cycle;
\draw[facestyle_r, color2,opacity = .5] (v7_r) -- (v8_r) -- (v9_r) -- (v7_r) -- cycle;

\draw[facestyle_r] (c1) -- (c2) -- (c3) -- (c4) -- cycle;
\draw[facestyle_r] (c3) -- (c4) -- (c8) -- (c7) -- cycle;
\draw[facestyle_r] (c2) -- (c3) -- (c7) -- (c6) -- cycle;

\draw[color1!80] (v0_r) -- (v1_r) -- (v6_r) -- (v11_r);
\draw[color2!80] (v7_r) -- (v9_r);
\draw[color0!80] (v2_r) -- (v3_r) -- (v4_r); 

\draw[thick,->] (0,0,0) -- (3,0,0) node[anchor=south]{$x_{01}$};
\draw[thick,->] (0,0,0) -- (0,3,0) node[anchor=north west]{$x_{02}$};
\draw[thick,->] (0,0,0) -- (0,0,3) node[anchor=north]{$x_{12}$};

\end{tikzpicture}
    \caption{
        The LP-relaxation of the tension polytope of our running example is cut with the three hyperplanes $S_0, S_1$ and $S_2$.
        The intersections are the embeddings of the three polytropes.
        }
    \label{fig:XLP} 
\end{figure}
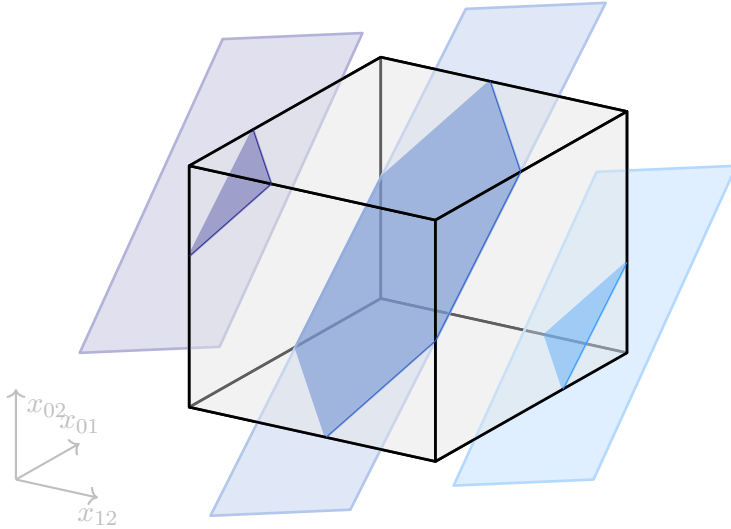

\begin{lemma}\label{lem:slicingTensionPolytope}
    For the affine subspace $S_z = \im B^\top + T p$ for some $p \in \Z^{A(G)}$ such that $z = \Gamma p,$ we have that
    \begin{equation}
        X_\text{LP} \cap S_z = \im \mathfrak{m}_z \subseteq X.
    \end{equation}
\end{lemma}
This property follows directly from the definition of the tension polytopes $X$ and $X_\text{LP}$. 
It allows us to view $X$ from a different perspective, as it is the convex hull of the pairwise disjoint embeddings of $\mathfrak{R}(z)$.
The affine subspaces $S_z$ are simply the same hyperplane translated at different points in space, cutting different sections of the cube $X_\text{LP}$.
We have:
\begin{corollary}\label{cor:embedding}
    The periodic tension polytope $X$ can be described by
    \begin{align}
        X   & = \conv \left\{ \im \mathfrak{m}_z \ \middle|\ \text{feasible~cycle~offset~} z\right\}.
    \end{align}
\end{corollary}
\Cref{fig:XLP} presents a visualization of these properties.
This observation motivates our interest in the number of affine subspaces $S_z$ or, equivalently, on the number of polytropes $\mathfrak{R}(z)$ in $\mathscr{T}$, which we will discuss in \Cref{sec:zonotopes}.

\subsection{Vertices of Timetabling Polytropes}

We now prove a result on the vertices of the polytropes, that mirrors a well-known property of the vertices of $X$, see \cite{nachtigallHabilitation}.
\begin{lemma}\label{lem:verticesAreSpanningSubgraphs}
    Each vertex of each polytrope in $\mathscr{T}$ corresponds to a unique spanning subgraph of $\overline{G}$, where $G$ is the digraph of the PESP instance.
\end{lemma}
\begin{proof}
    Let $\pi$ be a vertex of some polytrope $\mathbf{R}(p) \subset \mathscr{T}$.
    Being a point in $n-1$ dimensional space it must lay on at least $n-1$ distinct hyperplanes in the $\mathscr{H}$-description of $\mathbf{R}(p)$.
    As a weighted digraph polyhedron $R(p)$ is $\mathscr{H}$-described by lower and upper bound hyperplanes only, each one stemming from a different arc in $\overline{G}$.
    Let us then denote by $K$ the subgraph spanned by the arcs of $\overline{G}$ whose inequality is tight in $\pi$.
    Looking at the signed coordinate directions $\pm d_i$ for $i \in \{1, \ldots, n-1\}$, we have that for every $i$ at least one of the two directions $+d_i$ and $-d_i$ points outside of $\mathbf{R}(p)$ at $\pi$, for if that was not the case, then $\pi$ would not be a vertex to begin with.
    Since the spacial coordinates correspond to nodes in $V(\overline{G})$, this means that $x \coloneqq m_p(\pi)$ is either lower or upper bound tight for at least one arc incident in every node in $V(\overline{G})$.
    We can thereby construct $K$ as the subgraph induced by the collection of such arcs.
    Doing so, $K$ is a spanning subgraph and the tension $x$ is the tension that attains the lower and upper bounds in the arcs of $K$.
\end{proof}
Now an immediate consequence of \Cref{cor:embedding} together with how we constructed $K$ in the above proof is that vertices of polytropes are mapped via the maps $m_p$ to vertices of $X$.

To conclude we recall the known result of \cite{JoswigKulas}.
\begin{lemma}\label{lem:shortestPathTreesVertices}
    The tropical vertices of a polytrope are in bijection with shortest path trees rooted at each vertex of the digraph from which the polytrope arises as a weighted digraph polyhedron.
\end{lemma}

\begin{example}\label{ex:shortestPathTreesVertices}
    Considering the instance of \Cref{ex:running} with various periodic offsets $p$, we focus on the shortest paths rooted at $v_1$ in $(\overline{G}, \kappa(p))$.
    \Cref{fig:1rootedTrees} highlights them, and also specifies a representative of the timetable so obtained.
    The trees are coloured with the same colours of \Cref{fig:tiling}, and each timetable is indeed contained in the polytrope of corresponding offset and colour.
    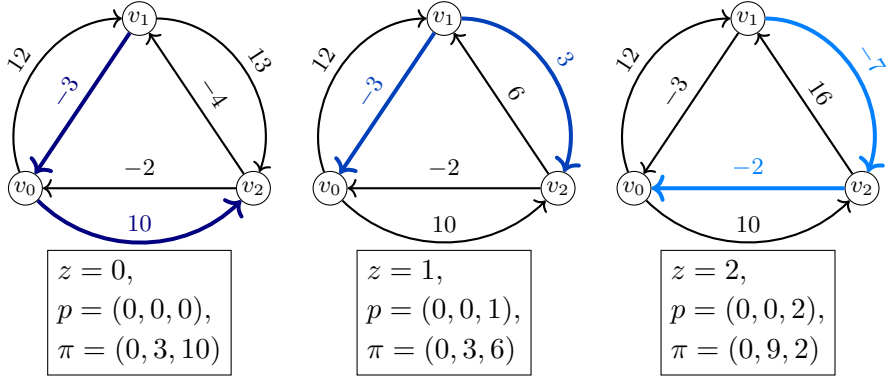
\begin{figure}
        \centering
        \definecolor{color0}{HTML}{000080}
\definecolor{color1}{HTML}{0040BB}
\definecolor{color2}{HTML}{0080FF}
\begin{tikzpicture}[scale=1.5]
    \tikzstyle{p} = [line width=1.5, ->]
    \tikzstyle{u} = [line width=1.5]
    \tikzstyle{q} = [line width=1.5, ->]
    \tikzstyle{r} = [line width=1.5, ->]
    \tikzstyle{v} = [draw, circle, inner sep=1, minimum width=12, font=\footnotesize]
    \tikzstyle{w} = [v]
    \tikzstyle{s} = [draw, rectangle, inner sep=1, minimum width=12, minimum height=12, font=\footnotesize]
    \tikzstyle{t} = [midway, font=\footnotesize, sloped]
    \tikzstyle{ta} = [t, above]
    \tikzstyle{tb} = [t, below]
    \tikzstyle{tr} = [t, right]
    \tikzstyle{tl} = [t, left]
    \tikzstyle{l} = [font=\scriptsize]
    \node[w] (A) at (0, 0) {$v_0$ };
    \node[v] (B) at (1, 1.5) {$v_1$};
    \node[w] (C) at (2, 0) { $v_2$};
    
    \draw[p, line width = 0.8] (A) edge[looseness=1,out=110,in=180] node[ta] {$12$} (B);
    \draw[p, line width = 0.8] (B) edge[looseness=1,out=0,in=70] node[ta] {$13$} (C);
    \draw[p, color0] (A) edge[looseness=1,out=-45,in=-135] node[ta] {$10$} (C);	
    
    \draw[p, color0] (B) -- node[ta] {$-3$} (A);
    \draw[p, line width = 0.8] (C) -- node[ta,black] {$-4$} (B);
    \draw[p, line width = 0.8] (C) -- node[ta,black] {$-2$} (A);
    
    \node[rectangle,draw,align=left] at (1,-1.1) {$z=0,$\\ $p = (0,0,0),$\\ $\pi = (0,3,10)$};
\end{tikzpicture}
%
\begin{tikzpicture}[scale=1.5]
    \tikzstyle{p} = [line width=1.5, ->]
    \tikzstyle{u} = [line width=1.5]
    \tikzstyle{q} = [line width=1.5, ->]
    \tikzstyle{r} = [line width=1.5, ->]
    \tikzstyle{v} = [draw, circle, inner sep=1, minimum width=12, font=\footnotesize]
    \tikzstyle{w} = [v]
    \tikzstyle{s} = [draw, rectangle, inner sep=1, minimum width=12, minimum height=12, font=\footnotesize]
    \tikzstyle{t} = [midway, font=\footnotesize, sloped]
    \tikzstyle{ta} = [t, above]
    \tikzstyle{tb} = [t, below]
    \tikzstyle{tr} = [t, right]
    \tikzstyle{tl} = [t, left]
    \tikzstyle{l} = [font=\scriptsize]
    \node[w] (A) at (0, 0) {$v_0$ };
    \node[v] (B) at (1, 1.5) {$v_1$};
    \node[w] (C) at (2, 0) { $v_2$};
    
    \draw[p, line width = 0.8] (A) edge[looseness=1,out=110,in=180] node[ta] {$12$} (B);
    \draw[p, color1] (B) edge[looseness=1,out=0,in=70] node[ta] {$3$} (C);
    \draw[p, line width = 0.8] (A) edge[looseness=1,out=-45,in=-135] node[ta] {$10$} (C);	
    
    \draw[p, color1] (B) -- node[ta] {$-3$} (A);
    \draw[p, line width = 0.8] (C) -- node[ta,black] {$6$} (B);
    \draw[p, line width = 0.8] (C) -- node[ta,black] {$-2$} (A);
    
    \node[rectangle,draw,align=left] at (1,-1.1) {$z=1,$\\ $p = (0,0,1),$\\ $\pi = (0,3,6)$};
\end{tikzpicture}
%
\begin{tikzpicture}[scale=1.5]
    \tikzstyle{p} = [line width=1.5, ->]
    \tikzstyle{u} = [line width=1.5]
    \tikzstyle{q} = [line width=1.5, ->]
    \tikzstyle{r} = [line width=1.5, ->]
    \tikzstyle{v} = [draw, circle, inner sep=1, minimum width=12, font=\footnotesize]
    \tikzstyle{w} = [v]
    \tikzstyle{s} = [draw, rectangle, inner sep=1, minimum width=12, minimum height=12, font=\footnotesize]
    \tikzstyle{t} = [midway, font=\footnotesize, sloped]
    \tikzstyle{ta} = [t, above]
    \tikzstyle{tb} = [t, below]
    \tikzstyle{tr} = [t, right]
    \tikzstyle{tl} = [t, left]
    \tikzstyle{l} = [font=\scriptsize]
    \node[w] (A) at (0, 0) {$v_0$ };
    \node[v] (B) at (1, 1.5) {$v_1$};
    \node[w] (C) at (2, 0) { $v_2$};
    
    \draw[p, line width = 0.8] (A) edge[looseness=1,out=110,in=180] node[ta] {$12$} (B);
    \draw[p, color2] (B) edge[looseness=1,out=0,in=70] node[ta] {$-7$} (C);
    \draw[p, line width = 0.8] (A) edge[looseness=1,out=-45,in=-135] node[ta] {$10$} (C);	
    
    \draw[p, line width = 0.8] (B) -- node[ta] {$-3$} (A);
    \draw[p, line width = 0.8] (C) -- node[ta,black] {$16$} (B);
    \draw[p, color2] (C) -- node[ta] {$-2$} (A);
    
    \node[rectangle,draw,align=left] at (1,-1.1) {$z=2,$\\ $p = (0,0,2),$\\ $\pi = (0,9,2)$};
     
\end{tikzpicture}
        \caption{Spanning tree structures rooted at $v_1$ in $\overline{G}$ for the three polytropes.}
        \label{fig:1rootedTrees}
    \end{figure}
\end{example}

With \Cref{lem:verticesAreSpanningSubgraphs} we additionally proved that non-tropical vertices of the polytrope, so-called \emph{pseudo-vertices} in \cite{JoswigKulas}, correspond to spanning trees that are not shortest path trees.
Notice that in reality we always spoke of spanning subgraph and not just spanning trees.
This is because of possible degeneracy, since even in full dimensional polytropes we may have that certain vertices are contained in more than just $d-1$ facets.

\subsection{The Tropical Neighbourhood}
\label{subsec:neighbourhood}

Having understood that the PESP timetable space is a finite collection of disjoint polytropes $\mathfrak{R}(z)$, each arising from a weighted digraph $(\overline{G}, \kappa(p))$ for some $p$ such that $\Gamma p = z$, the immediate next questions are to understand how these polytropes are actually positioned with respect to each other inside $\mathscr{T}$, and subsequently how the objective function of PESP behaves over said collection.
Recall that the objective function is $w \cdot x$, for some weights $w$ provided with the instance (note that these have nothing to do with the weights $\kappa(p)$ of the weighted digraphs), and as we know any tension $x$ depends heavily on the offset variables, as in 
\begin{equation}
    x = -B^\top \pi + Tp, \quad \text{~or similarly~} \quad \Gamma x = Tz.
\end{equation}
This implies that evaluating the objective function will be linear within a polytrope, where the offset is clearly constant, but non-linear in general, as changing polytrope and thereby changing integral offset vector will disrupt any hope of linearity.

\begin{figure}
    \centering
    	\begin{tikzpicture}[scale=.5]
	\definecolor{color0}{HTML}{000080}
	\definecolor{color1}{HTML}{0040BB}
	\definecolor{color2}{HTML}{0080FF}
	\tikzstyle{l} = [draw, line width=0.8, grau]
	\tikzstyle{v} = [circle, draw, black, line width=0.8, fill=white, inner sep=1.5]
	\tikzstyle{f} = []
	\tikzstyle{t} = [white, font=\tiny,  ]
	\tikzstyle{g} = [orange, dashed, line width=0.8]
	\draw[->] (0, -1) -- (15, -1);
	\draw[->] (0, -1) -- (0, 14);

	\foreach \i in {0,1,...,7,8, 9,10,...,14} {
		\draw (\i,-1.2) -- (\i,-0.8) node[below,yshift=-1ex,font=\tiny] {$\i$};
	};
	\foreach \i in {-1,0,1,2,3,4,...,13} {
		\draw (-0.2,\i) -- (0.2,\i) node[left,xshift=-1ex,font=\tiny] {$\i$};
	};
	
	\node[anchor=north, yshift=-0.5ex, font=\tiny] at (15,-1) {$\pi_1$};
	\node[anchor=east, font=\tiny] at (0,14) {$\pi_2$};
	\begin{pgfonlayer}{bg}
	\clip (0,-1) rectangle (15,14);
	\fill[f, fill=color0] (3,7) -- (3,12) -- (8,12) -- cycle;
	\node[t, anchor=west] at (3.25, 9) {0\,0\,0};
	\fill[f, fill=color0] (13,7) -- (13,12) -- (18,12) -- cycle;
	\node[t] at (14, 9) {-1\,0\,1};
	\fill[f, fill=color0] (3,-3) -- (3,2) -- (8,2) -- cycle;
	\node[t, anchor=west] at (3.25, -.5) {0\,1\,1};
	\fill[f, fill=color0] (13,-3) -- (13,2) -- (18,2) -- cycle;
	\node[t] at (14, -.5) {-1\,1\,2};
	
	\fill[f, fill=color1] (3,2) -- (3,7) -- (8,12) -- (13,12) -- (13,7) -- (8,2) -- cycle;
	\node[t] at (8, 6.5) {0\,0\,1};
	\fill[f, fill=color1] (13,2) -- (13,7) -- (18,12) -- (23,12) -- (23,7) -- (18,2) -- cycle;
	\node[t] at (14, 3) {-1\,0\,2};
	\fill[f, fill=color1] (-7,2) -- (-7,7) -- (-2,12) -- (3,12) -- (3,7) -- (-2,2) -- cycle;
	\node[t] at (1, 9) {1\,0\,0};
	
	\fill[f, fill=color1] (3,12) -- (3,17) -- (8,22) -- (13,22) -- (13,17) -- (8,12) -- cycle;
	\node[t] at (6, 13) {0\,-1\,0};
	\fill[f, fill=color1] (13,12) -- (13,16) -- (17,20) -- (22,20) -- (22,16) -- (18,12) -- cycle;
	\node[t] at (14, 13) {-1\,-1\,1};
	
	\fill[f, fill=color1] (3,-8) -- (3,-3) -- (8,2) -- (13,2) -- (13,-3) -- (8,-8) -- cycle;
	\node[t] at (10, -0.5) {0\,1\,2};
	\fill[f, fill=color1] (-7,-8) -- (-7,-3) -- (-2,2) -- (3,2) -- (3,-3) -- (-2,-8) -- cycle;
	\node[t] at (1, -0.5) {1\,1\,1};
	
	\fill[f, fill=color2] (8,2) -- (13,7) -- (13,2) -- cycle;
	\node[t, anchor=east] at (12, 3) {0\,0\,2};
	\fill[f, fill=color2] (-2,2) -- (3,7) -- (3,2) -- cycle;
	\node[t, anchor=east] at (2, 3) {1\,0\,1};
	
	\fill[f, fill=color2] (8,12) -- (13,17) -- (13,12) -- cycle;
	\node[t, anchor=east] at (12, 13) {0\,-1\,1};
	\fill[f, fill=color2] (-2,12) -- (3,17) -- (3,12) -- cycle;
	\node[t, anchor=east] at (2, 13) {1\,-1\,0};

	\end{pgfonlayer}

	\draw[g] (0,0) -- (10,0);%
	\draw[g] (0,10) -- (10,10);
	\draw[g] (0,0) -- (0,10);
	\draw[g] (10,0) -- (10,10);
\end{tikzpicture}
    \caption{
        Decomposition of the timetable space $\Pi'/\R \mathbf{1}$ as a covering of $\R^{V(G)}/\R \mathbf{1}$ corresponding to the limit case of \Cref{ex:running}.
        }
    \label{fig:tilingPrime}
\end{figure}
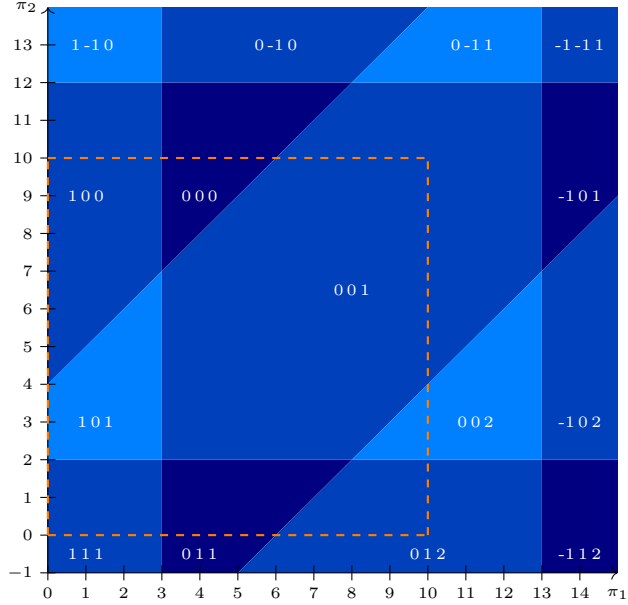

A natural first step is to see how a change in the offset vector shifts the focus from a certain polytrope in the torus to another. 
Unsurprisingly, we will see that the smaller this change is, the smaller the shift will be. 
To this end, we consider a specific and somewhat artificial PESP instance: 
\begin{definition}
    For a PESP instance $(G, T, \l, u, w)$, we define its \emph{limit instance} as 
    \begin{equation}
        (G, T, \l, w)_{\infty} \coloneqq (G, T, \l, \l + T, w).
    \end{equation}
    We denote $\mathbf{R}'(p)\subseteq \mathscr T,  R'(p)\subseteq \Pi'$ as the polytropes in the torus and the timetable space respectively, where $\Pi'$ corresponds to the timetable space for the limit case.
\end{definition}
Such construction is actually excluded a priori in our hypothesis, since we had $u - \l < T$, which is a reasonable assumption in practice.
Nonetheless, we consider $\Pi'$ anyway, since it helps us in presenting and understanding certain properties of the polytrope arrangement which we will then easily adapt to the more general and practical case.
\begin{lemma}
    For the limit case $(G, T, \l, w)_{\infty}$, the collection of polytropes $\mathbf{R}'(p)$ and $R'(p)$, for all $p \in \Z^{A(G)}$, cover the entire torus $\mathscr{T}$ and space $\R^{V(G)} = \Pi'$, respectively. 
\end{lemma}
This is true because one can choose any $\pi \in \R^{V(G)}$ and due to the span of $T$ for each $a \in A(G)$ we can easily find an integer $p_a$  to satisfy the bounds. 
In \Cref{fig:tilingPrime} we can see a depiction of the tiling of $\Pi'/\R\mathbf{1}$ corresponding to the limit case of the PESP instance of \Cref{ex:running}.
Comparing it to the original tiling in \Cref{fig:tiling}, we notice that certain "infeasible bands" that encircled the polytropes there now have disappeared in the tiling of $(G,T,\l,w)_\infty$.
Said bands are clearly of the form $u_{ij} - Tp_{ij} < \pi_j - \pi_i < \l_{ij} - T(p_{ij} - 1)$ for each arc $(i,j) \in A(G)$ and $p_{ij} \in \Z$, and their Manhattan-width is the $T$-complement of the span of the arc $(i,j)$, namely $T - (u_{ij} - \l_{ij})$.
In the case of $\Pi'$ these bands do not exist, since by construction we now set $u_{ij} - \l_{ij} = T$.
By that we understand that the decomposition of $\Pi'$ by polytropes is the one induced by a hyperplane arrangement, where all hyperplanes are of the form $\pi_j - \pi_i = \l_{ij} + Tp_{ij}$, for all possible integer values of $p_{ij}$.

\begin{figure}[htbp]
    \centering
    	\begin{tikzpicture}[scale=2.5]
						\tikzstyle{p} = [line width=1.5, ->]
						\tikzstyle{u} = [line width=1.5]
						\tikzstyle{q} = [line width=1.5, ->]
						\tikzstyle{r} = [line width=1.5, ->]
						\tikzstyle{v} = [draw, circle, inner sep=1, minimum width=12, font=\footnotesize]
						\tikzstyle{w} = [v]
						\tikzstyle{s} = [draw, rectangle, inner sep=1, minimum width=12, minimum height=12, font=\footnotesize]
						\tikzstyle{t} = [midway, font=\footnotesize, sloped]
						\tikzstyle{ta} = [t, above]
						\tikzstyle{tb} = [t, below]
						\tikzstyle{tr} = [t, right]
						\tikzstyle{tl} = [t, left]
						\tikzstyle{l} = [font=\scriptsize]
						\node[w] (I) at (0, 0) {$i$};
						\node[w] (J) at (2, 0) {$j$};
						
						\draw[p] (I) edge[looseness=1,out=15,in=165] node[ta] {$\l+T-p_{ij}T-kT$} (J);
						\draw[p] (I) edge[looseness=1,out=55,in=125] node[ta] {$\l+T-p_{ij}T$} (J);
						\draw[p] (J) edge[looseness=1,out=195,in=-15] node[tb] {$p_{ij}T+kT - \l$} (I);
						\draw[p] (J) edge[looseness=1,out=235,in=-55] node[tb] {$p_{ij}T - \l$} (I);

					\end{tikzpicture}
    \caption{
        The subgraph of $\overline{G'}$ induced by the sole vertices $i$ and $j$, prior to the resolution of the tropical sum of $(\overline{G'}, \kappa(p) \oplus \kappa(p + ke_{ij}))$.
        Note that we applied the substitution $u = \l + T$.
        }
    \label{fig:emptyIntersectionLemma}
\end{figure}
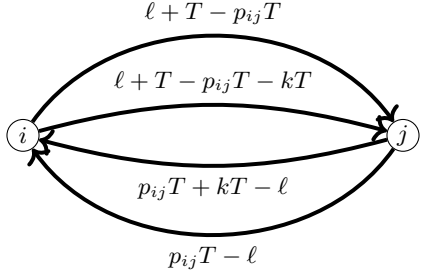

\begin{lemma}
    Let $p \in \Z^{A(G)}$ be an offset vector, $k$ an integer such that $k \notin \{-1,0,1\}$, and $e_{ij}$ the canonical basis vector of the arc $(i,j)$ . 
    Then we have that $R'(p) \cap R'(p + ke_{ij}) = \varnothing$.
\end{lemma}
\begin{proof}
    Considering $R'(p)$ and $R'(p + ke_{ij})$ as weighted digraph polyhedra, let us denote their weighted digraphs as $(\overline{G'}, \kappa(p))$ and $(\overline{G'}, \kappa(p + ke_{ij}))$ respectively.
    Now, making quick use of \cite[Lemma~6]{JoswigLoho}, we have that the intersection $R'(p) \cap R'(p + ke_{ij})$ is a weighted digraph polyhedron, whose weighted digraph is $(\overline{G'}, \kappa(p) \oplus \kappa(p + ke_{ij}))$.
    The only place in $\overline{G'}$ where the latter tropical sum is not obvious by idempotency is precisely between nodes $i$ and $j$, where we have a situation akin to that of \Cref{fig:emptyIntersectionLemma}.
    The first and last of the four arcs are the ones with weight $\kappa(p)$, whereas the second and third have weight $\kappa(p + ke_{ij})$.
    If $k > 1$ then the minimum $\kappa_{ij}$ is attained in the second arc and the minimum $\kappa_{ji}$ is attained in the fourth arc, and so these two arcs are kept.
    Now the weight of the remaining 2-cycle is $\l+T-p_{ij}T-kT+p_{ij}T-\l = T(1-k) < 0$, which is negative.
    If instead $k < -1$ we similarly keep the first and third arc, and again we have a negative weight cycle.
    In either case, $R'(p) \cap R'(p + ke_{ij}) = \varnothing$.
\end{proof}
If instead of having big jumps in the offset we only move by 1, we have the following.
\begin{lemma}\label{lem:neighboursDistanceOne}
    Let $p \in \Z^{A(G)}$ be an offset vector, $k$ an integer such that $k \in \{-1,1\}$, and $e_{ij}$ the canonical basis vector of the arc $(i,j)$ . 
    Then we have that $R'(p) \cap R'(p + ke_{ij}) = F$, where $F$ is the face of either polytrope in which one of the inequalities given by the arc $(i,j)$ is tight.
\end{lemma}
\begin{proof}
    Everything repeats as in the previous proof, except that now the weight of the remaining 2-cycle is $\l+T-p_{ij}T-kT+p_{ij}T-\l = T(1-k) = 0$ if $k = 1$ and $\l+T-p_{ij}T+p_{ij}T+kT-\l = T(1+k) = 0$ if $k = -1$.
    In either case, we did not create an immediate infeasibility in $i-j-i$, and our construction is identical to that in \cite[Lemma~7]{JoswigLoho}.
    By that lemma we have restricted to the face in which one of the inequalities given by the arc $(i,j)$ is tight, in particular the lower bound inequality if $k = 1$ and the upper bound inequality if $k = -1$.
    Note that $F$ can be the whole of $R'(p)$ as well as the empty face.
\end{proof}
To conclude we can better discuss the dimension of $F$, and gain some insight.
We operate with the same hypothesis of \Cref{lem:neighboursDistanceOne}, and briefly assume $R'(p)$ to be full-dimensional.
If $F$ is a facet of $R'(p)$, then it is also a facet of $R'(p + ke_{ij})$, which must then be full-dimensional as well.
If that was not the case we would have $R'(p + k e_{ij}) = F$, of dimension $n-1$, and this could only happen if there was a 0-weight cycle in $(\overline{G'}, \kappa(p))$ of length 2, which is incompatible with our assumption of $G$ having no antiparallel arcs.
When $F$ instead is not a facet of $R'(p)$, then $R'(p + k e_{ij}) = F$.
If that was not the case there would be a point in the interior of $R'(p + k e_{ij})$ that is separated from $R'(p)$ by the relevant $(i,j)$-inequality, which in turn would have to be facet-defining by Farkas' Lemma.

We can now revert back to the standard $\Pi$ and $\mathbf{R}(p)$ case.
Since for all arc spans $u_{ij} - \l_{ij} < T$ we are back to the case where all polytropes are disjoint.
By the preceding lemmas though, we now understand that only when an inequality arising from arc $(i,j)$ is facet-tight in $\mathbf{R}(p)$ can we hope for $\mathbf{R}(p \pm e_{ij})$ to be non-empty, since otherwise we have proven that we would be restricting $\mathbf{R}(p \pm e_{ij})$ to being a face of $\mathbf{R}(p)$, and now the non-zero Manhattan-width of the corresponding $(i,j)$-band would render $\mathbf{R}(p \pm e_{ij})$ empty.
Of course $\mathbf{R}(p \pm e_{ij})$ can still be empty anyway, when the complement of the span of arc $(i,j)$ is too large, meaning the infeasible band is too wide.

Looking at the tiling of the limit case, it is now reasonable to consider two polytropes to be neighbours if they share a common facet.
This leads us to extending this notion to the general case.
\begin{definition}
    Two non-empty polytropes $\mathbf{R}(p_1),\mathbf{R}(p_2) \subseteq \mathscr{T}$ are \emph{neighbours} if and only if the polytropes $\mathbf{R}'(p_1)$ and $\mathbf{R}'(p_2)$ of the limit case intersect in a facet.
    The \emph{neighbourhood of $\mathbf{R}(p)$} is the set of its non-empty neighbours, or
    \begin{equation}
        N_p = \left\{ \mathbf{R}(q) \neq \varnothing \ \middle| \ \exists (i,j) \in A(G) \colon p-q = \pm e_{ij} \right\}.
    \end{equation}
\end{definition}
In our example, we therefore consider $\mathbf{R}(0,0,1)$ to be a neighbour of both $\mathbf{R}(0,0,0)$ and $\mathbf{R}(0,0,2)$ while the two triangular polytropes are not neighbours, since $\mathbf{R}'(0,0,0)$ and $\mathbf{R}'(0,0,2)$ share no facet, as is evidenced in \Cref{fig:tilingPrime}.

The considerations above allow us to  describe neighbouring polytropes in a simpler way: we can obtain the neighbourhood relation with the periodic offsets instead of taking a detour via the limit case.
\begin{lemma}
    Two non-empty polytropes $\mathbf{R}(p_1)$ and $\mathbf{R}(p_2)$ are neighbours whenever there exist representatives $p_1'$ and $p_2'$, meaning offset vectors such that $p_1-p_1', p_2-p_2' \in \ker\Gamma$, such that $p_1' - p_2' = \pm e_{ij}$ for some arc $(i,j) \in A(G)$.
\end{lemma}
By this we also recognise as neighbours polytropes $\mathfrak{R}(z_1)$ and $\mathfrak{R}(z_2)$ whenever the difference $z_1 - z_2$ is, up to sign, a column of $\Gamma$.

\begin{example}\label{ex:neighbourhood-graph-construction}
    If we were now to try to visualise the entire tiling of $\mathscr{T}$ and its spatial configuration, a possible idea is to make an undirected graph $N$, whose nodes are feasible cyclic offsets $z \in \Z^\mathscr{B}$, and two nodes $z_1$ and $z_2$ are adjacent if and only if the corresponding polytropes $\mathfrak{R}(z_1)$ and $\mathfrak{R}(z_2)$ are neighbours.
    Notice that such a procedure is equivalent to making the nodes be periodic offsets $p \in \Z^A$, adjacent to their neighbours, and then taking the quotient graph with respect to our usual equivalence relation
    \begin{equation}
        p_1 \equiv p_2 \iff z_1 \coloneqq \Gamma p_1 = \Gamma p_2 \eqqcolon z_2.
    \end{equation}
    For our running example such a graph would be that of \Cref{fig:fullBorough}.
    \begin{figure}[htbp]
        \centering
        	\begin{tikzpicture}[scale=2.5]
    	\definecolor{color0}{HTML}{000080}
        \definecolor{color1}{HTML}{0040BB}
        \definecolor{color2}{HTML}{0080FF}
		\tikzstyle{p} = [line width=1.5, ->]
		\tikzstyle{u} = [line width=1.5]
		\tikzstyle{q} = [line width=1.5, ->]
		\tikzstyle{r} = [line width=1.5, ->]
		\tikzstyle{v} = [draw, circle, inner sep=1, minimum width=12, font=\footnotesize]
		\tikzstyle{w} = [v]
		\tikzstyle{s} = [draw, rectangle, inner sep=1, minimum width=12, minimum height=12, font=\footnotesize]
		\tikzstyle{t} = [midway, font=\footnotesize, sloped]
		\tikzstyle{ta} = [t, above]
		\tikzstyle{tb} = [t, below]
		\tikzstyle{tr} = [t, right]
		\tikzstyle{tl} = [t, left]
		\tikzstyle{l} = [font=\scriptsize]
		
		\node[w, fill = color0, text= white] (0) at (0,0) {$0$};
		\node[w, fill = color1, text= white] (1) at (1,0) {$1$};
		\node[w, fill = color2, text= white] (2) at (2,0) {$2$};
		
		\draw[u] (0) edge (1);
		\draw[u] (1) edge (2);
	
	\end{tikzpicture}
        \caption{
            The $N$ graph of the tiling of the initial example from \Cref{fig:ean}.
            The colours correspond to the coloured tiles in \Cref{fig:tiling}.
            }
        \label{fig:fullBorough}
    \end{figure}
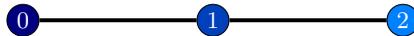
\end{example}

The construction presented in \Cref{ex:neighbourhood-graph-construction} serves as the basis of a local improvement procedure.
Given any known feasible solution $(\pi,x,p)$, we focus on the offset $p$, meaning a node in $N$.
Now the neighbourhood of $p$ in $N$ can be explored, and the PESP instance optimized for each offset in $p' \in N_{p}$.
This amounts to solving standard linear programs, one per each polytrope $R(p')$.
If a solution $(\pi',x',p')$ is found, the search continues in $N_{p'}$.
We can name this search procedure the tropical neighbourhood search, or \tns.
Such a strategy is in spirit quite similar to the state of the art approach, the modulo network simplex (also \mns), see \cite{ModuloNetworkSimplex}.
The major problem of \mns~is that it sometimes reaches local optima, meaning solutions $(\pi,x,p)$ from which no improvement can be found, even though $(\pi,x,p)$ is not the optimal solution overall.
Examining example instances and comparing the search space of \mns~with that of \tns~we have found some interesting differences.
First of all, there are adjacency relations in the \mns~search space that do not appear in the \tns~search space, and vice versa.
This means that while \mns~can, in theory, explore solutions starting from $(\pi,x,p)$ that are not in a neighbour of $\mathbf{R}(p)$, also \tns~can find among the neighbours $\mathbf{R}(p')$ for $p' \in N_p$ solutions that \mns~is not be able to consider.
We also found instances where the \mns~search space presents local optima whereas the \tns~search space has none.
\begin{example}\label{ex:mango}
    For a certain instance with 8 nodes, 10 arcs, and $\mu = 3$, we computed the neighbourhood graph for the \tns~and \mns~heuristics, respectively shown on the left and on the right of the figure.
    For \tns~each node corresponds to a polytrope in the polytropal decomposition, and each arc adjuncts two polytropes whose offsets differ by 1 in a single coordinate.
    For \mns~each node corresponds to a spanning tree of the PESP instance, and each arc adjuncts two spanning trees that differ only by removing one arc and adding one other.
    The redder a node, the better its objective value, with purple nodes being global optima.
    Squares, indeed here only present for \mns, are local optima.
    \begin{figure}
        \centering
        \includegraphics[width=0.45\textwidth,height=\textheight,keepaspectratio]{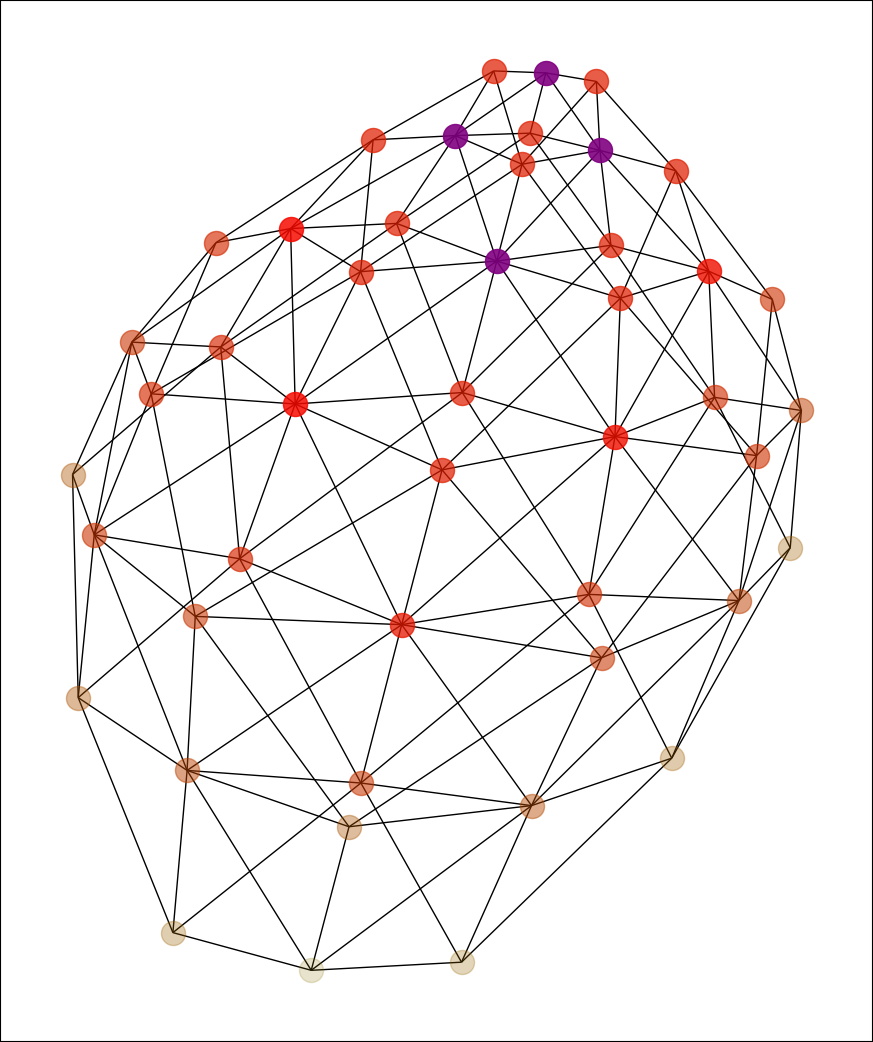}
        \includegraphics[width=0.45\textwidth,height=\textheight,keepaspectratio]{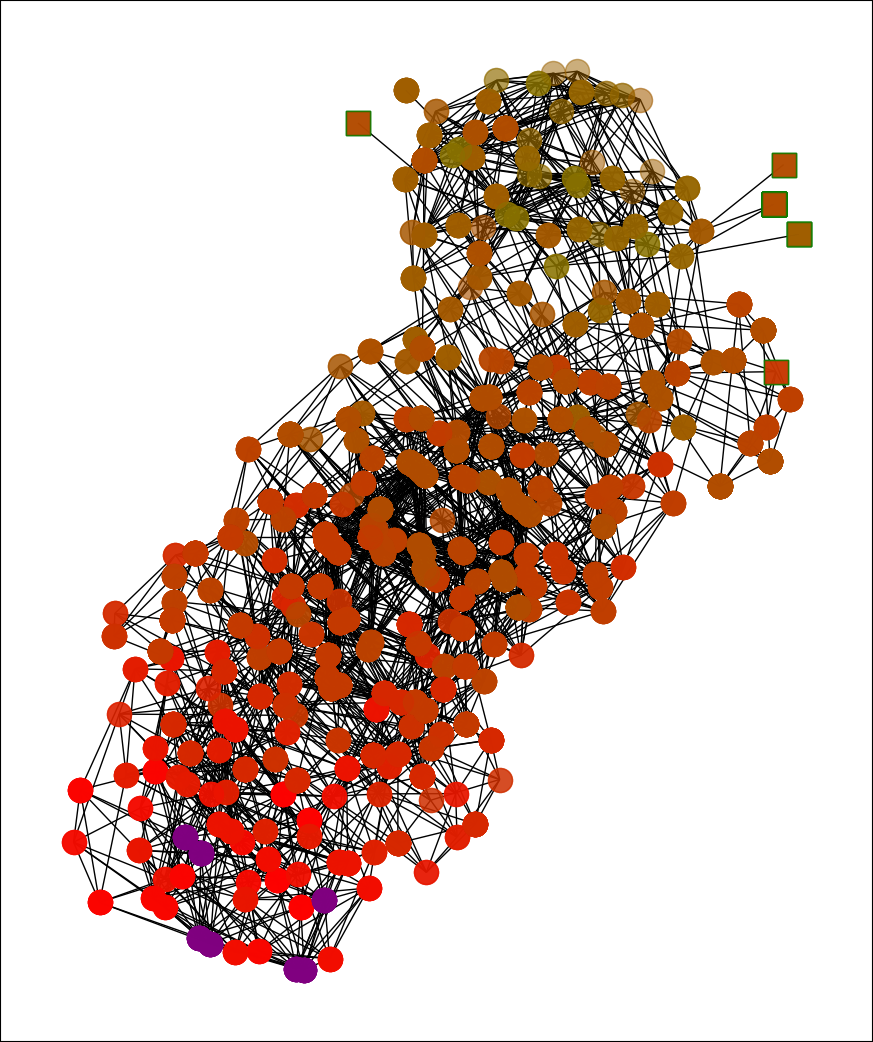}
        \caption{
            The neighbourhoods of the two heuristics compared.
            }
        \label{fig:mango}
    \end{figure}
\end{example}
In principle, this means that \tns~can be used twofold: either as a local improvement heuristic on its own, or as a tool to escape local minima in existing approaches. 
This brief analysis allows us to be hopeful for \tns~to be beneficial in practice.

\section{Cycle Offset Zonotopes}
\label{sec:zonotopes}

A \emph{zonotope} is the image of a hypercube under an affine map. 
In formulae, a zonotope in $\mathbb R^\mu$ can hence be written as 
\begin{equation}
    Z(A, b) \coloneqq \{ Ax \mid \mathbf{0} \leq x \leq \mathbf{1} \} + b
\end{equation}
for some matrix $A$ with $\mu$ rows and some translation vector $b \in \mathbb R^\mu$. 
We will abbreviate $Z(A) \coloneqq Z(A, \mathbf{0})$.

Recall that for a PESP instance $(G, T, \l, u, w)$ and an integral cycle basis $\mathscr{B}$ of $G$ with cycle matrix $\Gamma$, we defined the \emph{cycle offset zonotope} in \Cref{def:zonotope} as
        \begin{equation}
        Z \coloneqq \left\{ z \in \R^\mathscr{B} \ \middle|\ \exists x \in \R^{A(G)} \colon \Gamma x = T z, \ \l \leq x \leq u \right \}.
    \end{equation}
Let $\Gamma'$ be the matrix arising from $\Gamma$ by scaling each column with $\frac{u_a - \l_a}{T}$. Then
\begin{equation}
    Z = \left\{ \frac{\Gamma x}{T} \ \middle|\ \l \leq x \leq u \right\} \\
    =  \left\{ \Gamma' x \ \middle|\ \mathbf{0} \leq x \leq \mathbf{1} \right\} + \frac{\Gamma \l}{T} = Z\left(\Gamma', \frac{\Gamma \l}{T} \right),
\end{equation}
so that $Z$ is indeed a zonotope.

We assume for this chapter that $u_a - \l_a > 0$ holds for all $a \in A(G)$. 
This can always be achieved by contracting arcs with $\l_a = u_a$ \cite{LiebchenBook}. 
As a consequence, the rank of $\Gamma'$ is $\mu$, and the dimension of $Z$ is the cyclomatic number $\mu = |\mathscr{B}|$.

\subsection{Zonotopal Tilings, Spanning Trees and the Number of Polytropes}

It is a well-known property of zonotopes that they can be tiled by subzonotopes, see \cite[\S7.5]{crt-beck}.
\begin{definition}
    A \emph{zonotopal tiling} of a zonotope $Z(A)$ is a polyhedral subdivision of $Z(A)$ such that the cells are translates of zonotopes $Z(A_S)$, where $A_S$ is the submatrix arising from $A$ by deleting the columns indexed by $S$.
    We call each maximal cell of such subdivision a \emph{tile}.
    The tiling is called \emph{fine} if all tiles are parallelotopes, i.e., the generators of each tile are linearly independent.
\end{definition}
Every zonotope has a fine tiling \cite[Corollary~7.5.10]{crt-beck}. 
If $A$ has $\mu$ rows and full row rank, then the tiles of a fine tiling are given as certain translates of $Z(A_S)$, where $A_S$ ranges over all invertible $(\mu \times \mu)$-submatrices of $A$.

\begin{lemma}\label{lem:tile-spanning-tree}
    Let $Z$ be a cycle offset zonotope of a PESP instance on the digraph $G$. 
    Fix any fine zonotopal tiling. 
    Then there is a one-to-one correspondence between spanning subgraphs of $G$ with $|A(G)| - k$ arcs and the $k$-dimensional cells of the tiling such that the arc set $S$ of a spanning subgraph maps to a translate of $Z(\Gamma'_S)$.
\end{lemma}
\begin{proof}
    A $(\mu \times \mu)$-submatrix of $\Gamma$ and hence $\Gamma'$ is invertible if and only if its columns correspond to the $\mu$ co-tree arcs of a spanning tree of $G$ \cite[Lemma~3.1]{kavitha_cycle_2009}. 
    This establishes the correspondence for $k = \mu$. 
    If $S$ is a general spanning subgraph with $|A(G)| - k$ arcs, then it contains a spanning tree, so that the corresponding submatrix $\Gamma'_S$ is a $(\mu \times k)$-submatrix of an invertible $(\mu \times \mu)$-submatrix.
    This implies that $\Gamma'_S$ has rank $k$ and hence $Z(\Gamma'_S)$ is $k$-dimensional. 
    Conversely, if some $(\mu \times k)$-submatrix $\Gamma'_S$ has rank $k$, then we can add further columns from $\Gamma'$ so that $\Gamma'_S$ is a submatrix of an invertible $(\mu \times \mu)$-matrix, which corresponds to a spanning tree. 
    It follows that $S$ contains a spanning tree and is therefore spanning.
\end{proof}

\begin{example}
\label{ex:smallTiling}
    The cycle offset zonotope of \Cref{ex:running} is rather simple, as shown in \Cref{fig:runningZonotope}. 
    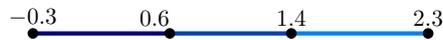
\begin{figure}[htbp]
        \centering
        \begin{tikzpicture}[scale = 2]  
    
 \definecolor{color0}{HTML}{000080}
	\definecolor{color1}{HTML}{0040BB}
	\definecolor{color2}{HTML}{0080FF}
    \tikzstyle{vertex} = [text=black, inner sep=0.5pt, align=left, draw=none, scale = 0.8] 

    \draw[color0, ultra thick] (-0.3, 0) -- (0.6, 0);
     \draw[color1, ultra thick] (0.6, 0) -- (1.4, 0);
      \draw[color2, ultra thick] (1.4, 0) -- (2.3, 0);
    \fill[black] (-0.3, 0) circle (1 pt);
    \fill[black] (0.6, 0) circle (1 pt);
    \fill[black] (1.4, 0) circle (1 pt);
    \fill[black] (2.3, 0) circle (1 pt);
   
    \node[vertex, above] at (-0.3, 0.04) {$-0.3$};
    \node[vertex, above] at (0.5, 0.04) {$0.6$};
    \node[vertex, above] at (1.4, 0.04) {$1.4$};
    \node[vertex, above] at (2.3, 0.04) {$2.3$};
    
\end{tikzpicture}
        \caption{The cycle offset zonotope of \Cref{ex:running}, with a fine tiling.}
        \label{fig:runningZonotope}
    \end{figure}
    Given that the graph has only one cycle the dimension is indeed just 1.
    The translation value is easily computed to be 
    \begin{equation}
        \frac{\Gamma \l}{T} = \frac{\begin{pmatrix} 1 & -1 & 1 \end{pmatrix}\begin{pmatrix} 3 \\ 2 \\ 4 \end{pmatrix}}{10} = 0.5 ,
    \end{equation}
    and then the tiles are given by invertible submatrices of $\Gamma'$.
    Since 
    \begin{equation}
        \Gamma' = \begin{pmatrix} 0.9 & -0.8 & 0.9 \end{pmatrix}
    \end{equation}
    there are clearly three possible such $(1 \times 1)$-submatrices, corresponding to the three entries taken individually.
    In \Cref{fig:runningZonotope} there are three line segments delimited by adjacent points, and each one corresponds to a tile.
    From the translation point $0.5$ we in fact have a backward tile of length $-0.8$, corresponding to the spanning tree $v_0$--$v_1$--$v_2$, and two forward tiles both of length $0.9$, corresponding to the spanning trees $v_0$--$v_2$--$v_1$ and $v_1$--$v_0$--$v_2$.
    In particular, this is a fine tiling of the zonotope.
    If instead of the two forward tiles we considered the tile induced by the subgraph $v_0$--$v_2$ we would have a single forward tile, generated by both $0.9$ generators together, going from $0.5$ to $2.3$.
\end{example}

The cycle offset zonotope $Z$ shares the property that tiles correspond to spanning trees with graphical zonotopes. 
However, our notion is somewhat dual: The lower-dimensional cells of the zonotopal tiling are given by spanning subgraphs in our case, whereas the cells of a fine zonotopal tiling of a graphical zonotope are given by forests. 

Having established this correspondence, we can obtain a formula for the volume of a cycle offset zonotope $Z$.
We denote the set of spanning trees of a graph $G$ by $\mathscr S(G)$.
We will often identify a spanning tree $S \in \mathscr S(G)$ with its set of arcs.
\begin{corollary}\label{cor:volume}
    The volume of a cycle offset zonotope $Z$ of a PESP instance $(G, T, \l, u, w)$ is
    \begin{equation}
        \vol(Z) = \sum_{S \in \mathscr S(G)} \prod_{a \in A(G) \setminus A(S)} \frac{u_a - \l_a}{T}.
    \end{equation}
    In particular, the volume does not depend on the choice of the cycle basis.
\end{corollary}
\begin{proof}
    By \Cref{lem:tile-spanning-tree}, we can compute the volume $\vol(Z)$ as the sum of the volumes $\vol(Z(\Gamma'_S))$, for $S \in \mathscr S(G)$. As any $Z(\Gamma'_S)$ is a parallelotope,
    \begin{equation}
        \vol(Z(\Gamma'_S)) = |\det(\Gamma'_S)| = \prod_{a \in A(G) \setminus S} \frac{u_a - \l_a}{T} \cdot |\det(\Gamma_S)| = \frac{u_a - \l_a}{T},
    \end{equation}
    where we used $|\det(\Gamma_S)| = 1$ because any invertible $(\mu \times \mu)$-submatrix of the cycle matrix $\Gamma$ has determinant $\pm 1$.
\end{proof}

\begin{theorem}\label{thm:lattice-spanning}
    Let $Z$ be a cycle offset zonotope of a PESP instance on $G$. 
    Then $Z$ has at most as many integer points as there are spanning trees of $G$.
\end{theorem}
\begin{proof}
    Let $A$ be an invertible $(\mu \times \mu)$-submatrix of the cycle matrix $\Gamma$. 
    We claim that $\mathbf{0}$ is the only integer point in the set $Q \coloneqq \{Ay \mid -\mathbf{1} < y < \mathbf{1} \}$. 
    Indeed, for any such integer point $q \in Q$, there is a unique solution $y$ to $Ay = q$. 
    By Cramer's rule, $y$ is integral, as the determinant of $A$ is $\pm 1$. 
    But then $y = \mathbf{0}$ and hence $q = \mathbf{0}$.
    
    Now consider a fine tiling of $Z$ and let $Z(A') + b$ be a tile. 
    That is, $A'$ is an invertible $(\mu \times \mu)$-submatrix of $\Gamma$, and $b$ is some translation vector, possibly non-integral. 
    Suppose that $z, z' \in Z(A') + b$ are both integral. 
    Then there are $x, x'$ with $\mathbf{0} \leq x, x' \leq \mathbf{1}$ such that $A'x + b = z$ and $A' x' + b = z'$. 
    Thus $A'(x - x')$ is integral. 
    Let $y$ denote the vector with entries $(x_a - x_a')(u_a - \l_a)/T$, for $a \in A(G)$. 
    If $A$ is the submatrix of $\Gamma$ corresponding to $A'$, then $Ay = A'(x - x')$ and $-\mathbf{1} < y < \mathbf{1}$ because $u_a - \l_a < T$. 
    We conclude that $y = 0$ and hence $x = x'$. 
    This means that any tile in a fine tiling of $Z$ contains at most one lattice point. 
    In particular, using \Cref{lem:tile-spanning-tree}, the number of integer points is bounded by $|\mathscr S(G)|$.
\end{proof}

\begin{example}
\label{ex:favex-tiled}
    \begin{figure}
        \centering
        \begin{tikzpicture}[scale=2.5]
					\tikzstyle{p} = [line width=1, ->]
					\tikzstyle{v} = [draw, circle, inner sep=1, minimum width=12, font=\footnotesize]
					\tikzstyle{t} = [midway, font=\footnotesize]
					\tikzstyle{ta} = [t, above]
					\tikzstyle{tb} = [t, below]
					\tikzstyle{tr} = [t, right]
					\tikzstyle{tl} = [t, left]
					
					\node[v] (A) at (0, 0) {$v_0$};
					\node[v] (B) at (0, 1) {$v_1$};
					\node[v] (C) at (1, 1) {$v_2$};
					\node[v] (D) at (1, 0) {$v_3$};

					\draw[p] (B) -- node[ta, black, sloped] {$[3,12]$} (A);
					\draw[p] (B) -- node[ta, black, sloped] {$[3,12]$} (C);
					\draw[p] (D) -- node[tb, black, sloped] {$[3,12]$} (C);	
					\draw[p] (D) -- node[ta, black, sloped] {$[3,12]$} (A);
					\draw[p] (A) edge[bend left] node[tb, black, sloped] {$[6,15]$} (B);
					\draw[p] (C) edge[bend left] node[tb, black, sloped] {$[4,13]$} (D);
					
\end{tikzpicture}
        \caption{A PESP instance with period $T = 10$.}
        \label{fig:favex}
    \end{figure}
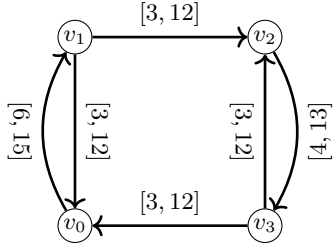
    Consider the PESP instance graph in \Cref{fig:favex}.
    Contrary to our hypothesis it contains antiparallel arcs, which we initially wanted to avoid for the sake of simplicity of certain dimension arguments in \Cref{subsec:neighbourhood}.
    In general, these directed cycles of length two do not hinder neither PESP nor the heart of our considerations.
    In particular, we can construct the cycle offset zonotope of this instance and tile it.
    As integral cycle basis, we use the cycles induced by the three regions of the planar embedding as in \Cref{fig:favex}.
    The tiled zonotope is shown in \Cref{fig:favexZonotopeTiled}.
    \begin{figure}
        \centering
        \input{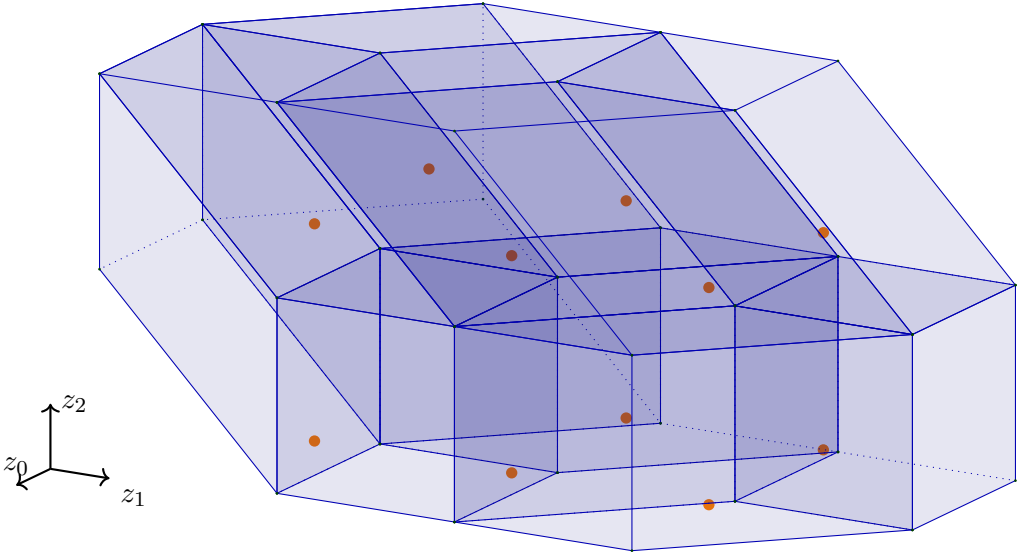}
        \caption{The tiled zonotope of the PESP instance of \Cref{fig:favex}, with integer points marked in orange.}
        \label{fig:favexZonotopeTiled}
    \end{figure}
    There are 12 tiles in the picture, but only 11 integer points inside the zonotope, marked in orange. 
    This is coherent with \Cref{lem:tile-spanning-tree} and \Cref{thm:lattice-spanning}.
    In fact the graph of \Cref{fig:favex} has exactly 12 spanning trees, and so 12 tiles are found.
    Moreover each tile has the potential of containing at most 1 integer point, and in fact we have 11 such points.
    In truth here there is one integer point that lays on a common facet of two tiles, thereby being contained in both.
    Regardless, it is in general still possible for a tile to contain no integer point altogether.
\end{example}

\subsection{A Duality Between Torus Polytropes and Zonotope Tiles}

Let $(G, T, \l, u, w)$ be a PESP instance. 
In \Cref{sec:tropical}, we viewed the space of feasible periodic timetables modulo symmetries as subset of a torus $\mathscr T$. 
Fixing an integral cycle basis $\mathscr B$, the timetable space on the torus decomposes into the union of pairwise disjoint polytropes $\mathfrak R(z)$ for all $z \in \Z^{\mathscr B}$. 
\begin{lemma}
    \label{lem:polytrope-to-lattice}
    The map $z \mapsto \mathfrak R(z)$ is a bijection between integer points of the cycle offset zonotope $Z$ and non-empty polytropes in $\mathscr T$.
\end{lemma}
\begin{proof}
    Let $z$ be an integer point of $Z$. 
    Then, by definition of $Z$, there is a vector $x$ with $\l \leq x \leq u$ such that $\Gamma x = T z$. 
    By \Cref{lem:slicingTensionPolytope}, $x \in \im \mathfrak{m}_z$, and so $\mathfrak R(z) \neq \varnothing$. 
    Conversely, if $\mathfrak R(z)$ is a non-empty polytrope in $\mathscr T$, then we find $x \in \im \mathfrak{m}_z$, and $\frac{\Gamma x}{T} = z$ is an integer point of $Z$.
\end{proof}

\begin{figure}[htbp]
    \centering
    \begin{tikzpicture}[scale = 2]  
     \definecolor{color0}{HTML}{000080}
	\definecolor{color1}{HTML}{0040BB}
	\definecolor{color2}{HTML}{0080FF}
    \tikzstyle{vertex} = [text=black, inner sep=0.5pt, align=left, draw=none, scale = 0.8] 

    \draw[color0, ultra thick] (-0.3, 0) -- (0.6, 0);
     \draw[color1, ultra thick] (0.6, 0) -- (1.4, 0);
      \draw[color2, ultra thick] (1.4, 0) -- (2.3, 0);
    \fill[black] (-0.3, 0) circle (1 pt);
    \fill[black] (0.6, 0) circle (1 pt);
    \fill[black] (1.4, 0) circle (1 pt);
    \fill[black] (2.3, 0) circle (1 pt);
    \fill[orange] (0, 0) circle (1.2 pt);
    \fill[orange] (1, 0) circle (1.2 pt);
    \fill[orange] (2, 0) circle (1.2 pt);
    \node[vertex, above] at (-0.3, 0.04) {$-0.3$};
    \node[vertex, above] at (0.5, 0.04) {$0.6$};
    \node[vertex, above] at (1.4, 0.04) {$1.4$};
    \node[vertex, above] at (2.3, 0.04) {$2.3$};
    
\end{tikzpicture}
    \caption{The tiled cycle offset zonotope of \Cref{ex:running}, with integer points marked in orange.}
    \label{fig:runningZonotopeWithPoints}
\end{figure}
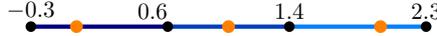

By \Cref{thm:lattice-spanning}, we immediately obtain:
\begin{theorem}
    The number of polytropes in the decomposition of the timetable space on $\mathscr T$ is at most $|\mathscr S(G) |$.
\end{theorem}
Looking at \Cref{fig:runningZonotopeWithPoints}, we clearly see that for the instance of \Cref{ex:running} there are exactly three integral points in the zonotope.
In fact, there are also three tiles, three polytropes as seen in \Cref{fig:torus}, and three possible spanning trees of the graph shown in \Cref{fig:ean}.

In view of \Cref{lem:polytrope-to-lattice}, the ``maximal'' objects of the polytropal decomposition of $\mathscr T$ hence correspond to certain $0$-dimensional objects of $Z$. 
We will investigate now how polytrope vertices as $0$-dimensional objects relate to the tiles of a zonotopal tiling as top-dimensional objects of $Z$. 
We understand these two relations as a kind of duality, but we do not expect that this extends to all dimensions in between, as the polytropes have dimension at most $|V(G)| - 1$, whereas the zonotope tiles have dimension $\mu$.

We will start with investigating the faces of the fractional periodic tension polytope $X_\text{LP} = \{ x \in \R^{A(G)} \mid \l \leq x \leq u \}$.
\begin{definition}
    A \emph{structure} $(S, L, U)$ in $G$ is a triplet of subsets of $A(G)$ such that $S = L \cup U$ and $L \cap U = \varnothing$.
    If $S$ induces a spanning subgraph of $G$, we say $(S, L, U)$ is a \emph{spanning structure}.
    If the subgraph is also a tree, we say $(S, L, U)$ is a \emph{spanning tree structure}.
\end{definition}
We will use structures as a combinatorial notation for the faces of $X_\text{LP}$, which is justified by the following observation:
\begin{lemma}
    The map
    \begin{align}
        (S, L, U) \mapsto \mathcal{F}_{L,U} \coloneqq \left\{ x \in X_\text{LP} \ \middle | \ x_a = \l_a \ \forall a \in L, x_a = u_a \ \forall a \in U \right\}
    \end{align}
    is an inclusion-reversing bijection between structures in $G$ and faces of $X_\text{LP}$. 
    A structure with $|S| = k$ corresponds to a face of dimension $|A(G)| - k$.
\end{lemma}

The cycle offset zonotope $Z$ is by definition the image of $X_\text{LP}$ under $\frac{1}{T}\Gamma$.
As the faces $\mathscr F_{L,U}$ are (scaled) hyperrectangles themselves, $\frac{1}{T}\Gamma$ maps each  $\mathscr F_{L,U}$ to a subzonotope of $Z$. 
A straightforward computation shows the following:
\begin{lemma}
    \label{lem:gamma-on-cube-face}
    Let $(S, L, U)$ be a structure. Then
    \begin{align}
        \left\{\frac{\Gamma x}{T} \ \middle| \ x \in \mathscr F_{L,U} \right\} =
      Z\left (\Gamma'_S, \frac{\Gamma v}{T} \right),
    \end{align}
    where $v_a \coloneqq u_a$ for $a \in U$ and $v_a \coloneqq \l_a$ otherwise. 
\end{lemma}

\begin{lemma}
    \label{lem:cells-from-structures}
    Let $C$ be a cell of a fine zonotopal tiling of $Z$ such that $C$ is a translate of $Z(\Gamma'_S)$ for some $S \subseteq A(G)$. 
    Then there exist $L$ and $U$ such that $(S, L, U)$ is a structure and $C$ is the image of $\mathscr F_{L, U}$ under $\frac{1}{T}\Gamma$.
\end{lemma}
\begin{proof}
    In view of \Cref{lem:gamma-on-cube-face}, this is only a question of translation vectors. 
    By \cite{barcelo_zonotopal_1994}, keeping in mind that our zonotopes are not centrally symmetric w.r.t.\ the origin, we can turn a zonotopal tiling in our sense to a so-called strong zonotopal tiling such that the translation vectors are exactly of the form as in \Cref{lem:gamma-on-cube-face}.
\end{proof}

We are now ready to formulate a duality result between tiles of a zonotopal tiling and polytrope vertices.
\begin{theorem}
    \label{thm:tile-to-vertex}
    Fix an arbitrary fine zonotopal tiling of $Z$. 
    Then any tile is the image of $\mathscr F_{L,U}$ under $\frac{1}{T}\Gamma$ for a spanning tree structure $(S, L, U)$.
    If a tile contains an integer point $z \in Z$, then the polytrope $\mathfrak R(z)$ has a vertex defined by $(S, L, U)$.
\end{theorem}
\begin{proof}
    The first statement follows from \Cref{lem:tile-spanning-tree} and \Cref{lem:cells-from-structures}. 
    If a tile corresponding to $(S, L, U)$ contains a lattice point $z$, then we find $x \in \mathscr F_{L, U}$ with $\Gamma x = Tz$. 
    But then $x$ is a vertex of $X$ and by \Cref{lem:verticesAreSpanningSubgraphs} also a vertex of $\mathfrak{R}(z)$.
\end{proof}

\begin{example}
    Consider the middle tile in \Cref{fig:runningZonotopeWithPoints}.
    As we will see in \Cref{fig:P-polyhedron} it arises from the spanning tree structure that takes the lower bound in the arc $(v_0, v_1)$ and the upper bound in the arc $(v_1, v_2)$. 
    As the tile clearly contains an integer point we know it is also related to a feasible polytrope in the polytropal decomposition.
    The associated polytrope is, in fact, the hexagon in \Cref{fig:tiling}, and here we see that the polytrope contains the vertex $(3,6)$, which corresponds to the same spanning tree. 
    In particular, the spanning tree is a shortest path tree rooted at $v_1$ in \Cref{fig:goverline}, and the vertex of the polytrope is a tropical vertex.
\end{example}

Now, to complement \Cref{thm:lattice-spanning}, we have the following statement.
\begin{corollary}
    If all the arc bounds of the PESP instance are free, i.e. $u_a - \l_a = T - 1$ for every arc $a$, and there is no degeneracy, meaning that all spanning structures for which a feasible solution exists are spanning tree structures, then the bound expressed in \Cref{thm:lattice-spanning} is tight.
\end{corollary}
\begin{proof}
    It is well known that if all arcs are free, then for any spanning tree structure there exists a feasible solution \cite{preProGoerigk}.
    This implies that given an instance with only free arcs, then any tile of any tiling of its cycle offset zonotope $Z$ must contain an integer point.
    Now, if an integer point in $Z$ is shared by two distinct tiles, this implies by \Cref{lem:cells-from-structures} that the two corresponding distinct spanning tree structures are simultaneously tight for the same associated solution, thereby giving the kind of degeneracy we excluded in our hypothesis.
    We have therefore shown that the tiling of a free instance that has no degeneracy must have each tile containing a point in its interior, concluding the proof.
\end{proof}

\subsection{Constructing Zonotopal Tilings}
\label{subsec:construction}

It is natural to ask for a converse of \Cref{thm:tile-to-vertex}: Suppose that for each non-empty polytrope $\mathfrak R(z)$ we pick a vertex corresponding to a spanning tree structure $(S, L, U)$ and map $\mathscr F_{L, U}$ to the cycle offset zonotope $Z$ via $\frac{1}{T}\Gamma$. 
When we choose pairwise distinct spanning trees per polytrope, do these $\mu$-dimensional zonotopes form the tiles of a zonotopal tiling of $Z$? 
We will describe a construction that produces a zonotopal tiling from certain spanning tree structures.

\begin{definition}
\label{def:P_i}
    Let $(G, T, \l, u, w)$ be a PESP instance and let $i \in V(G)$. 
    Denote by $B$ the incidence matrix of $G$. 
    We define the unbounded polyhedron
    \begin{equation}
        P_i \coloneqq X_\text{LP} + \pos(\widehat{B^\top}_i),
    \end{equation}
    where $\widehat{B^\top}_i$ is $B^\top$ without the $i$-th column, $\pos(\widehat{B^\top}_i)$ indicates the positive cone of $\widehat{B^\top}_i$, and the summation is the Minkowski sum.
\end{definition}

To start it is practical to recall that $\operatorname{rank} B^\top = n-1$, where $n \coloneqq |V(G)|$. 
Moreover we assumed the graph $G$ to be connected and thereby every column $b_i$ of $B^\top$ is non-zero, and so we have that $-b_i \in \pos(\widehat{B^\top}_i)$ and $\widehat{B^\top}_i$ has full rank for every $i$. 
From this we conclude that $\pos(\widehat{B^\top}_i)$ is a $(n-1)$-dimensional pointed cone.

Let $\mathscr B$ be an integral cycle basis of $G$ with cycle matrix $\Gamma$, and let $Z$ be the cycle offset polytope. 
Observe that $\frac{1}{T}\Gamma$ maps $P_i$ to $Z$, due to $X_\text{LP} \subset P_i$ and \Cref{thm:cycle-periodicity}. 
It will turn out that the bounded faces of $P_i$ produce a fine zonotopal tiling of $Z$, the $\mu$-dimensional bounded faces giving the tiles. 
We need a few preparatory lemmata to prove this result.

\begin{lemma}
    Every bounded face of $P_i$ is a face of $X_\text{LP}$.
\end{lemma}
\begin{proof}
    Let $F$ be a bounded face of $P_i$. 
    Then there are $d \in \R^{A(G)}$ and $e \in \R$ such that $d^\top x \geq e$ is valid for $P_i$ and $F = \{x \in P_i \mid d^\top x = e\}$. 
    Since $X_\text{LP} \subseteq P_i$, $d^\top x \geq e$ is also valid for $X_\text{LP}$, so that $F' \coloneqq \{x \in X_\text{LP} \mid d^\top x = e\}$ is a face of $X_\text{LP}$. 
    Clearly $F' \subseteq F$. 
    Conversely, if $x \in F$, then there are $x' \in X_\text{LP}$ and $p \in \pos(\widehat{B^\top}_i)$ such that $x = x' + p$. 
    Then $d^\top x' + d^\top p= e$ and since $d^\top x' \geq e$, $d^\top p \leq 0$. 
    But then $d^\top x + \lambda d^\top p \leq e$ for every $\lambda \in \R_{\geq 0}$, and $d^\top x + \lambda d^\top p \geq e$ due to $x + \lambda p \in P_i$. 
    We conclude that $x + \lambda p \in F$ for all $\lambda \in \R_{\geq 0}$. 
    Since $F$ is bounded, we must have $p = 0$ and hence $x = x' \in F'$. 
    This shows $F \subseteq F'$.
\end{proof}

We can give a combinatorial description of the bounded faces of $P_i$:

\begin{lemma}
    \label{lem:cone-flow-condition}
    A face $\mathscr F_{L,U}$ of $X_\text{LP}$ is a bounded face of $P_i$ if and only if there is a vector $f \in \R_{\geq 0}^{A(G)}$ such that $f_a > 0$ for all $a \in L \cup U$, $f_a = 0$ otherwise, and
    \begin{equation}\label{eq:cone-flow-condition}
        \forall j \in V(G)\setminus\{i\} \colon \sum_{a \in \delta^+(j) \cap U} f_a - \sum_{a \in \delta^+(j) \cap L} f_a - \sum_{a \in \delta^-(j) \cap U} f_a + \sum_{a \in \delta^-(j) \cap L} f_a < 0.
    \end{equation}
\end{lemma}
\begin{proof}
    Let $(S, L, U)$ be a structure. 
    We observe first that the inequality 
    \begin{equation}
        \sum_{a \in U} f_a x_a - \sum_{a \in L} f_a x_a \leq \sum_{a \in U} f_a u_a - \sum_{a \in L} f_a \l_a
    \end{equation}
     is valid for $X_\text{LP}$, and 
    \begin{equation}
        \mathscr F_{L, U} = \left \{ x \in X_\text{LP} \ \middle| \ \sum_{a \in U} f_a x_a - \sum_{a \in L} f_a x_a = \sum_{a \in U} f_a u_a - \sum_{a \in L} f_a \l_a \right \}
    \end{equation}
    for any vector $f \in \R_{\geq 0}^{A(G)}$ with $f_a > 0$ for all $a \in S$ and $f_a = 0$ otherwise. 
    
    If $\mathscr F_{L, U}$ is a bounded face of $P_i$, then $f$ can be chosen in such a way that
    \begin{equation}\label{eq:pos-condition}
        \sum_{a \in U} f_a x_a - \sum_{a \in L} f_a x_a  < 0
    \end{equation}
    for all generators $x$ of the cone $\pos(\widehat{B^\top}_i)$.
    Conversely, if we find such a $f$ satisfying \eqref{eq:pos-condition}, then $\mathscr F_{L, U}$ is bounded in $P_i$. 
    
    The cone $\pos(\widehat{B^\top}_i)$ is generated by the columns of $B^\top$ except the $i$-th column. 
    By definition of the incidence matrix, \eqref{eq:pos-condition} is equivalent to \eqref{eq:cone-flow-condition}.
\end{proof}

Recall that we constructed the graph $\overline{G}$ in \Cref{sec:tropical}, where we added reverse copies of arcs to $G$. 
In this language, \Cref{lem:cone-flow-condition} states:

\begin{corollary}\label{cor:pos-flow-g-bar}
    A face $\mathscr F_{L,U}$ of $X_\text{LP}$ is a bounded face of $P_i$ if and only if there is a positive flow with negative balances at all vertices $v \in V(G) \setminus \{i\}$ in the subgraph $\overline{G}_{L, U}$ of $\overline{G}$ with $A(\overline{G}_{L, U}) \coloneqq \{(j, k) \mid (j, k) \in U\} \cup \{(k, j) \mid (j, k) \in L\}$ and $V(\overline{G}_{L, U}) \coloneqq V(\overline{G}) = V(G)$.
\end{corollary}

Let $\mathscr F_{L, U}$ be a bounded face of $P_i$. 
Then we find a flow as in the statement of \Cref{cor:pos-flow-g-bar} in the subgraph $\overline{G}_{L, U}$. 
Since the balances of every flow sum up to $0$ on each component of $\overline{G}_{L, U}$, and only the balance at node $i$ can be negative, we must have that $\overline{G}_{L, U}$ is connected and spanning. 
This has the following consequence:

\begin{corollary}\label{cor:bounded-spanning}
    Let $\mathscr F_{L, U}$ be a bounded face of $P_i$. 
    Then the corresponding structure $(S, L, U)$ is spanning. 
    In particular, the dimension of $\mathscr F_{L, U}$ is at most $\mu$.
\end{corollary}

\begin{theorem}\label{thm:arborescence}
    There is a bijection between the $\mu$-dimensional bounded faces of $P_i$ and arborescences in $\overline{G}$ rooted at $i$, given by $\mathscr F_{L, U} \mapsto \overline{G}_{L,U}$.
\end{theorem}
\begin{proof}
    Let $(S, L, U)$ be a structure such that $\overline{G}_{L,U}$ is an arborescence rooted at $i$. 
    Define $f_{(j, k)} \coloneqq n^{n - d_k}$, where $n = |V(\overline{G}_{L,U})| = |V(G)|$ and $d_j$ is the distance from $i$ to $j$ in $\overline{G}_{L,U}$, i.e., the number of edges of the unique $i$-$j$-path in $\overline{G}_{L, U}$. 
    Since every vertex $j \neq i$ has exactly one ingoing arc and at most $n-1$ outgoing arcs in $\overline{G}_{L,U}$,
    \begin{equation}
        \sum_{a \in \delta^+(j)} f_a -  \sum_{a \in \delta^-(j)} f_a = |\delta^+(j)| n^{n-d_j-1} - n^{n-d_j} <  n \cdot n^{n-d_j-1} - n^{n-d_j} = 0.
    \end{equation}
    Here there is a positive flow in $A(\overline{G}_{L,U})$ with negative balance at all vertices except $i$, and we conclude by \Cref{cor:pos-flow-g-bar} that $\mathscr F_{L, U}$ is a bounded face of $P_i$.
    
    For the converse, suppose that $\mathscr F_{L, U}$ is a bounded face of $P_i$ of dimension $\mu$. 
    Then \Cref{cor:pos-flow-g-bar} yields the existence of a positive flow $f$ in the spanning tree $\overline{G}_{L, U}$ with negative balances at all vertices except $i$. 
    For $j \in V(\overline{G}_{L, U})$, let $d_j$ again denote its distance from $i$ (the number of edges of the unique undirected $i$-$j$-path in $\overline{G}_{L,U}$). 
    Let $d_\text{max} \coloneqq \max\{d_j \mid j \in V(\overline{G}_{L, U})\}$. 
    For $j \neq i$, we denote by $\text{pred}(j)$ the unique vertex $k$ with $d_k = d_j - 1$ to which $j$ is adjacent.
    
    We prove by induction on $d_\text{max} - d_j$ that for all vertices $j \neq i$ the arc $(\text{pred}(j), j)$ is the unique ingoing arc at $j$ in $\overline{G}_{L, U}$. 
    This is sufficient to show that $\overline{G}_{L, U}$ is an arborescence. 
    In the case $d_j = d_\text{max}$, $j$ is a leaf of the tree $\overline{G}_{L, U}$. 
    Since the balance of $f$ is negative, the unique arc incident with $j$ must be ingoing at $j$. 
    If $0 < d_j < d_\text{max}$, then all arcs that connect $j$ with vertices of higher distance must be outgoing at $j$ by induction hypothesis. 
    Since the balance of $f$ is negative, the unique arc that connects $j$ with a vertex of lower distance must be ingoing, so that $(\text{pred}(j), j)$ is indeed the unique ingoing arc at $j$.
\end{proof}

\begin{lemma}\label{lem:parallelotope}
    For any $k$-dimensional bounded face $\mathscr F_{L, U}$ of $P_i$, the image $\frac{1}{T}\Gamma(\mathscr F_{L, U})$ is a $k$-dimensional parellolotope.
\end{lemma}
\begin{proof}
    Since $\mathscr F_{L, U}$ is a hyperrectangle, it is enough to show that the restriction of $\Gamma$ to $\mathscr F_{L, U}$ is injective.
    By \Cref{cor:bounded-spanning}, $\mathscr F_{L, U}$ corresponds to a spanning structure $(S, L, U)$.
    In particular, there is a spanning tree $S'$ contained in $S$. 
    If $\Gamma$ happens to be the cycle matrix of the fundamental cycles $\gamma_1, \dots, \gamma_\mu$ of $S'$, the injectivity is clear: For any $x \in \mathscr F_{L, U}$, the entries corresponding to arcs of $S'$ are the same, so that $x$ is determined by its entries for the co-tree arcs of $S'$. 
    But the $k$-th entry $\gamma_k^\top x$ of $\Gamma x$ is 
    \begin{equation}
        \sum_{a \in S': \gamma_a \neq 0} \gamma_{k, a} x_a + \sum_{a \notin S': \gamma_a \neq 0} \gamma_{k,a} x_a .
    \end{equation}
    The first sum is constant, and the second sum only has the summand $\gamma_{k, a'} x_{a'}$ for the unique co-tree arc $a'$ of $S'$ contained in the fundamental cycle $\gamma_k$. 
    Since there is one fundamental cycle for each co-tree arc, $\Gamma$ is injective on $\mathscr F_{L, U}$.
    
    If $\Gamma$ is another integral cycle basis, then there is a $\Z$-invertible matrix $M$ such that $M\Gamma$ is the cycle matrix of a fundamental cycle basis for $S'$. 
    Since $M\Gamma|_{\mathscr F_{L,U}}$ is injective by the above argument, so is $\Gamma|_{\mathscr F_{L,U}}$.
\end{proof}

\begin{theorem}
    The bounded faces of $P_i$ are mapped via $\frac{1}{T}\Gamma$ to the tiles of a fine zonotopal tiling of $Z$.
\end{theorem}
\begin{proof}
    Let $\mathscr F_i$ denote the union of all bounded faces of $P_i$. 
    As $X_\text{LP} \subseteq P_i$ and $\Gamma (\pos(\widehat{B^\top}_i)) = 0$ by \Cref{thm:cycle-periodicity}, $\frac{1}{T} \Gamma(\mathscr F_i) = \frac{1}{T} \Gamma(P_i) = Z$. 
    In particular,
    \begin{equation}\label{eq:vol-bounded}
    \begin{aligned}
        \vol(Z) &= \vol \left( \bigcup_{\substack{\mathscr F_{L,U} \subseteq \mathscr F_i \\ \dim\mathscr F_{L,U} = \mu}} \frac{1}{T}\Gamma(\mathscr F_{L,U}) \right) \\
        &\leq \sum_{\substack{\mathscr F_{L,U} \subseteq \mathscr F_i \\ \dim\mathscr F_{L,U} = \mu}} \vol \left(\frac{1}{T}\Gamma(\mathscr F_{L, U}) \right).
    \end{aligned}
    \end{equation}
    By \Cref{thm:arborescence}, the $\mu$-dimensional bounded faces of $P_i$ correspond one-to-one with arborescences in $\overline{G}_{L, U}$ rooted at $i$. 
    Any such arborescence yields a spanning tree structure $(S, L, U)$, and conversely, for every spanning tree $S$ of $G$, there is a unique spanning tree structure $(S, L, U)$ that turns $S$ into an arborescence on  $\overline{G}_{L, U}$. 
    With \Cref{lem:gamma-on-cube-face}, \eqref{eq:vol-bounded} implies
    \begin{equation}\label{eq:vol-spanning-trees}
        \vol(Z) \leq \sum_{S \in \mathscr{S}(G)} \vol(Z(\Gamma'_S)),
    \end{equation}
    and by the proof of \Cref{cor:volume}, \eqref{eq:vol-spanning-trees} and \eqref{eq:vol-bounded} are in fact equations. 
    In particular, the images of two distinct bounded $\mu$-dimensional faces $\mathscr F_{L, U}$ have an intersection of dimension strictly less than $\mu$. 
    It remains to remark that by \Cref{lem:parallelotope}, $\frac{1}{T}\Gamma(\mathscr F_{L,U})$ is always a parellolotope.
\end{proof}

We return to the decomposition of the torus $\mathscr T$ into polytropes $\mathfrak R(z)$. 
Any vertex of $\mathfrak R(z)$ is defined by a spanning tree structure by \Cref{lem:verticesAreSpanningSubgraphs}. 
Let $i \in V(G)$. 
Then the $i$-th tropical vertex is defined by the spanning tree structure corresponding to the shortest path arborescence rooted at $i$ in $\overline{G}_{L, U}$. 
For any $i \in V(G)$, we can hence assign to any polytrope $\mathfrak R(z)$ a spanning structure $(S^i_z, L^i_z, U^i_z)$ such that $\mathscr F_{L^i_z, U^i_z}$ is a bounded face of $P_i$. 
As a result, we obtain a partial converse to \Cref{thm:tile-to-vertex}:

\begin{theorem}
    Let $i \in V(G)$. 
    Then there is a fine zonotopal tiling of $Z$ such that any tile containing an integer point $z$ corresponds to a spanning tree structure defined by the $i$-th tropical vertex of $\mathfrak R(z)$.
\end{theorem}

\begin{example}
\label{fig:P-polyhedron}
    To illustrate the properties above, we again refer to our running example, \Cref{ex:running}. 
    Dropping the column of vertex $v_1$ of $B^T$ results in the unbounded polyhedron as displayed in \Cref{fig:cubePlusCone}.
    The coloured arrows depict the cone generators corresponding to the columns of vertices $v_0$ and $v_2$ in $B^T$. 
    The bounded faces are the three blue edges of the cube corresponding to the spanning tree structures $(S_i, L_i, U_i)$ for the pairs
    \begin{equation}
        \begin{aligned}
            L_0 = \{01\}, U_0 = \{02\},\\
            L_1 = \{01\}, U_1 = \{12\},\\
            L_2 = \{02\}, U_2 = \{12\},
        \end{aligned}
    \end{equation}
    respectively.
    Mapping each of those faces via $\frac{1}{T}\Gamma$ results in the three tiles as seen in \Cref{fig:runningZonotopeWithPoints}. 
    In this fine zonotopal tiling of $Z$, each tile contains an integer point, namely $z = 0$, $z = 1$ and $z = 2$ respectively. 
    As we saw with \Cref{lem:shortestPathTreesVertices}, the $i$-th tropical vertices in each $\mathfrak{R}(z)$ are identified with the shortest path arborescences rooted at $v_i$ in the corresponding weighted digraphs.
    In \Cref{ex:shortestPathTreesVertices} the $v_1$-rooted arborescences were discussed.
    The arcs used in each such arborescence identify the sets $L_i$ and $U_i$ that define the bounded face of the corresponding cycle offset.
    Again in \Cref{ex:shortestPathTreesVertices} three timetables arose.
    With the appropriate value of $z$ these three timetables are embedded via $\mathfrak{m}_z$ into $X_{\text{LP}}$, precisely to the three coloured bounded faces in \Cref{fig:cubePlusCone}.
    
    \begin{figure}
        \centering
        \begin{tikzpicture}[
    x = {(0.7cm,0.4cm)},
    y = {(0cm,1cm)},
    z = {(0.9cm,-0.2cm)},
    scale = 0.4,
    color = {lightgray}]

    \definecolor{darkorange}{HTML}{b37202}
    \definecolor{darkgreen}{HTML}{008500}
    \definecolor{color0}{HTML}{000080}
	\definecolor{color1}{HTML}{0040BB}
	\definecolor{color2}{HTML}{0080FF}
        \definecolor{pointcolor_r}{rgb}{ 1,0,0 }
        \tikzstyle{pointstyle_r} = [fill = pointcolor_r]
            \coordinate (c1) at (3,2,4);
            \coordinate (c2) at (3,2,13);
            \coordinate (c3) at (3,10, 13);
            \coordinate (c4) at (3, 10, 4); 
            \coordinate (c5) at (12, 2, 4);
            \coordinate (c6) at (12,2,13);
            \coordinate (c7) at (12,10, 13);
            \coordinate (c8) at (12, 10, 4); 
        \definecolor{edgecolor_r}{rgb}{ 0,0,0 }
        
        \definecolor{scolor}{rgb}{ 0.925,0.925,0.925 }
        
        \tikzstyle{facestyle_r} = [
            dotted,
            fill=blue!50!black,
            fill opacity=0.05, 
            line width = 0 pt, 
            line cap = round, 
            line join = round]
        \tikzstyle{backLines} = [
            dashed,
            fill = black, 
            fill opacity = 0, 
            line width = 0 pt, 
            line cap = round, 
            line join = round]
        \tikzstyle{frontLines} = [
            dashed,
            fill = black, 
            fill opacity = 0, 
            line width = 0.5 pt, 
            line cap = round, 
            line join = round]

    \draw[dashed] (c4) -- (c8);
    \draw[dashed] (c8) -- (c7);
    \draw[dashed] (c6) -- (c2);
    \draw[dashed] (c3) -- (c2);
    \draw[dashed] (c3) -- (c4);
    \draw[dashed] (c5) -- (c1);
    \draw[dashed] (c5) -- (c8);
    \draw[dashed] (c5) -- (c6);
    \draw[dashed, black] (c2) -- (c1);
    \draw[dashed, black] (c4) -- (c1);
    \draw[dashed, black] (c3) -- (c7);
    \draw[dashed, black] (c7) -- (c6);
    \draw[facestyle_r] (c2) -- (c3) -- (c4) -- ($(c4)+(0,-6,-6)$) -- ($(c2)+(0,-6,-6)$) -- cycle;
    \draw[facestyle_r] (c6) -- (c2) -- (c3) -- ($(c3)+(6,6,0)$) -- ($(c6)+(6,6,0)$) -- cycle;
    \draw[facestyle_r] (c6) -- (c2) -- ($(c2)+(0,-6,-6)$) -- ($(c6)+(0,-6,-6)$) -- cycle;
    \draw[facestyle_r] (c3) -- (c4) -- ($(c4)+(6,6,0)$) -- ($(c3)+(6,6,0)$) -- cycle;
    \draw[facestyle_r] (c4) -- ($(c4)+(6,6,0)$) -- ($(c4)+(0,-6,-6)$) -- cycle;
    
    \foreach \i in {3,4,6} {
        \draw[-, black, very thick, line cap = round] (c\i)-- ($(c\i)+(6,6,0)$);
        \draw[->, darkgreen, ultra thick, line cap = round] (c\i)-- ($(c\i)+(3,3,0)$);
        }
    \foreach \i in {2,4,6} {
        \draw[-, black, very thick, line cap = round] (c\i)-- ($(c\i)+(0,-6,-6)$);
        \draw[->, pink!50!red, ultra thick, line cap = round] (c\i) -- ($(c\i) + (0,-3,-3)$);
        }

    \draw[-, color0, line width = 2.5pt, line cap = round] (c4) -- (c3);
    \draw[-, color1, line width = 2.5pt, line cap = round] (c3) -- (c2);
    \draw[-, color2, line width = 2.5pt, line cap = round] (c2) -- (c6);
    \draw[thick,->] (0,0,0) -- (3,0,0) node[anchor = south]{$x_{01}$};
    \draw[thick,->] (0,0,0) -- (0,3,0) node[anchor = north west]{$x_{02}$};
    \draw[thick,->] (0,0,0) -- (0,0,3) node[anchor = north]{$x_{12}$};

\end{tikzpicture}
        \caption{The construction, according to \Cref{def:P_i}, of the object $P_1$ for the instance of \Cref{ex:running}.}
        \label{fig:cubePlusCone}
    \end{figure}
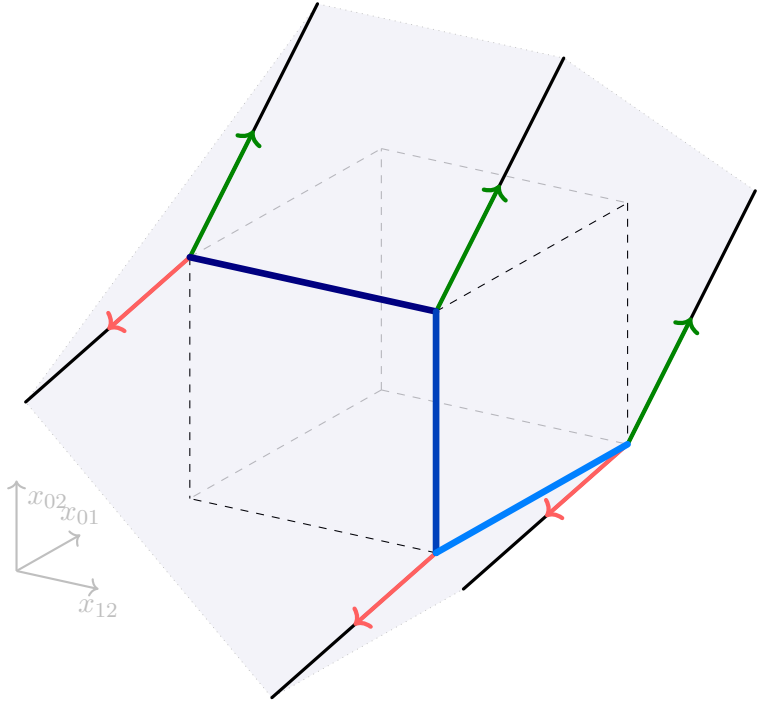
\end{example}

\subsection{Approximating the minimum width of a cycle basis}
\label{subsec:width}

For every choice of an integral cycle basis, the number of lattice points in the cycle offset zonotope $Z$ is the same by \Cref{lem:polytrope-to-lattice}. 
However, in order to solve the cycle-based MIP formulation for PESP \eqref{eq:pesp-formulation-cycle}, it has been observed that the choice of a cycle basis has an impact on the performance of branch-and-cut-based solvers.

\begin{definition}
    Let $(G, T, \l, u, w)$ be a PESP instance. 
    The \emph{width} of an integral cycle basis $\mathscr B$ of $G$ is defined as
    \begin{equation}
        W_{\mathscr{B}} \coloneqq \prod_{\gamma \in \mathscr B} \left( \left\lfloor \frac{\gamma_+^\top u - \gamma_-^\top \l}{T} \right\rfloor -\left\lceil \frac{\gamma_+^\top \l - \gamma_-^\top u}{T} \right\rceil +1 \right).
    \end{equation}
    Here, we decompose the oriented cycles $\gamma \in \mathscr B$ into their positive and negative parts $\gamma_+ \coloneqq \max(\mathbf{0}, \gamma)$ and $\gamma_- \coloneqq \max(\mathbf{0}, -\gamma)$, respectively.
\end{definition}

The rationale behind this notion is the following: For any feasible solution $(x, z) \in \R^{A(G)} \times \Z^{\mathscr B}$ to \eqref{eq:pesp-formulation-cycle} and for $\gamma \in \mathscr B$, it follows from $\l \leq x \leq u$ that
\begin{equation}\label{eq:odijk}
    \frac{\gamma_+^\top \l - \gamma_-^\top u}{T} \leq  z_\gamma = \frac{\gamma^\top x}{T} \leq \frac{\gamma_+^\top u - \gamma_-^\top \l}{T}.
\end{equation}
Since $z_\gamma$ is integer, we can round up the left-hand side and round down the right-hand side, so that there are at most
\begin{equation}\label{eq:rounding-argument}
    \left\lfloor \frac{\gamma_+^\top u - \gamma_-^\top \l}{T} \right\rfloor -\left\lceil \frac{\gamma_+^\top \l - \gamma_-^\top u}{T} \right\rceil +1
\end{equation}
possible integer values that the variable $z_\gamma$ can attain.

The width $W_{\mathscr B}$ is an upper bound on the number of lattice points in the cycle offset zonotope $Z$ for $\mathscr B$. 
Moreover, it is an upper bound on the number of leaves of the branch-and-bound-tree for the MIP \eqref{eq:pesp-formulation-cycle}, and this bound is sharp if the tree is fully explored. 
These considerations lead to the following optimization problem:
\begin{definition}
    Let $(G, T, \l, u, w)$ be a PESP instance. 
    The \emph{minimum width integral cycle basis problem} is to find an integral cycle basis $\mathscr B$ of minimum width $W_{\mathscr B}$.
\end{definition}
The complexity of this problem is open. 
It turns out that the minimum width integral cycle basis problem can be approximated by finding a minimum weight integral cycle basis, where the weight of an oriented cycle is 
\begin{equation}
    \sum_{a \in A(G) \colon \gamma_a \neq 0} (u_a - \l_a),
\end{equation}
see \cite{cyclePESPLiebchen}. 
Unfortunately, also this version has an unresolved complexity status, as finding a minimum weight fundamental cycle basis is NP-hard, whereas the problem is solvable in polynomial time on the larger class of undirected cycle bases \cite{horton_polynomial-time_1987}.

The number of lattice points in $Z$ is not only bounded by $W_{\mathscr B}$, but also by $|\mathscr{S}(G)|$, by \Cref{thm:lattice-spanning}. 
The aim of this section is to explore how far the number of spanning trees can be from the minimum width of an integral cycle basis.

Let $(G, T, \l, u, w)$ be a PESP instance, $\mathscr B$ an integral cycle basis of $G$, and let $Z$ be the corresponding cycle offset zonotope in $\R^{\mathscr B}$. 
By \eqref{eq:odijk}, we see that $Z$ is contained in the hyperrectangle
\begin{equation} 
    R \coloneqq \prod_{\gamma \in \mathscr B} \left[\frac{\gamma_+^\top \l - \gamma_-^\top u}{T}, \frac{\gamma_+^\top u - \gamma_-^\top \l}{T} \right].
\end{equation}
In fact, $R$ is the smallest hyperrectangle containing $Z$: For a specific $\gamma \in \mathscr B$, taking the vector $x \in X_\text{LP}$ with $x_a \coloneqq u_a$ if $\gamma_a > 0$ and $x_a \coloneqq \l_a$ otherwise produces a point $z \coloneqq \frac{\Gamma x}{T} \in Z$ with $z_\gamma = \frac{\gamma_+^\top u - \gamma_-^\top \l}{T}$, and the argument for the lower bound in the interval is analogous. 
With the rounding argument above \eqref{eq:rounding-argument}, we immediately find:
\begin{lemma}
    The width $W_{\mathscr B}$ equals the number of lattice points in $R$.
\end{lemma}

\begin{example}
    Consider the instance of \Cref{ex:favex-tiled}. 
    Recall that the number of lattice points of the cycle offset zonotope is $11$, as depicted in \Cref{fig:favexZonotopeTiled}.
    The width of the cycle basis induced by the three regions of the planar embedding is 
    \begin{equation}
        \begin{aligned}
        W_{\mathscr B} &= \left( \left \lfloor \frac{15 + 12}{10} \right \rfloor - \left \lceil \frac{6 + 3}{10} \right \rceil + 1 \right) \\
        &\quad \cdot 
        \left( \left \lfloor \frac{12 - 3 + 12 - 3}{10} \right \rfloor - \left \lceil \frac{3 - 12 + 3 - 12}{10} \right \rceil + 1 \right)\\
        &\quad \cdot
        \left( \left \lfloor \frac{12 + 13}{10} \right \rfloor - \left \lceil \frac{3 + 4}{10} \right \rceil + 1 \right)\\
        &= 2 \cdot 3 \cdot 2 = 12.
        \end{aligned}
    \end{equation}
    Indeed, $12$ is the number of lattice points in the hyperrectangle 
    \begin{equation}
        R = \left[ \frac{9}{10}, \frac{27}{10} \right]
        \times \left[ \frac{-18}{10}, \frac{18}{10} \right]
        \times \left[ \frac{7}{10}, \frac{25}{10} \right].
    \end{equation}
\end{example}

We hence want to compare the lattice points in $Z$ with the lattice points in $R$. We will do so by relating both with the volumes.
\begin{lemma}
\label{lem:volume-comparison}
    Let $d \in \R_{\geq 0}^{A(G)}$. 
    Then
    \begin{equation}
        \sum_{S \in \mathscr S(G)} \left (\prod_{a \in A(G) \setminus A(S)} d_a \right) \leq \prod_{\gamma \in \mathscr B} \left (\sum_{a \in A(G): \gamma_a \neq 0} d_a \right).
    \end{equation}
\end{lemma}
\begin{proof}
    Let $\Gamma$ be the cycle matrix of $\mathscr B$. 
    Proceeding as in the proof of \Cref{cor:volume}, the left-hand side is the volume of the zonotope $Z(\Gamma'')$, where $\Gamma''$ arises from $\Gamma$ by scaling the column corresponding to $a \in A(G)$ with $d_a$. 
    The hyperrectangle
    \begin{equation}
        \prod_{\gamma \in \mathscr B} \left[-\gamma_-^\top d, \gamma_+^\top d \right]
    \end{equation}
    contains $Z(\Gamma'')$, and its volume is
    \begin{equation}
        \prod_{\gamma \in \mathscr B} (\gamma_+ + \gamma_-)^\top d = \prod_{\gamma \in \mathscr B} \left (\sum_{a \in A(G): \gamma_a \neq 0} d_a \right),
    \end{equation}
    the right-hand side in the statement.
\end{proof}

A beautiful consequence of \Cref{lem:volume-comparison} for $d = \mathbf{1}$ is:
\begin{theorem}
    For any directed graph $G$ and any integral cycle basis $\mathscr B$ of $G$, the number of spanning trees in $G$ is at most $\prod_{\gamma \in \mathscr B} |\{a \in A(G) \mid \gamma_a \neq 0\}|$, i.e., the product of the lengths of the oriented cycles in $\mathscr B$.
\end{theorem}

Returning to PESP, we find:
\begin{theorem}
\label{thm:width-approx}
    Suppose that $W_{\mathscr B} \geq 1$. 
    Then
    \begin{equation}
        |\mathscr S(G)| \cdot  \left( \frac{\varepsilon}{T} \right)^\mu  \leq \vol(Z) \leq \prod_{\gamma \in \mathscr B} s_\gamma \leq W_{\mathscr B} \cdot \prod_{\gamma \in \mathscr B} \frac{s_\gamma}{\max\{\lfloor s_\gamma \rfloor, 1\}} < W_{\mathscr B} \cdot 2^\mu ,
    \end{equation}
    where $\varepsilon \coloneqq \min\{u_a - \l_a \mid a \in A(G)\}$, and
    \begin{equation}
        s_\gamma \coloneqq \sum_{a \in A(G): \gamma_a \neq 0} \frac{u_a - \l_a}{T}.
    \end{equation}
\end{theorem}
\begin{proof}
    By \Cref{cor:volume}, 
    \begin{equation}
        \vol(Z) = \sum_{S \in \mathscr S(G)} \left( \prod_{a \in A(G) \setminus A(S)} \frac{u_a - \l_a}{T} \right) \geq \sum_{S \in \mathscr S(G)} \left( \frac{\varepsilon}{T} \right)^\mu = |\mathscr S(G)| \cdot   \left( \frac{\varepsilon}{T} \right)^\mu .
    \end{equation}
    On the other hand, by \Cref{lem:volume-comparison}, 
    \begin{equation}
        \vol(Z) \leq \vol(R) = \prod_{\gamma \in \mathscr B} s_\gamma.
    \end{equation}
    Since $s_\gamma \geq 0$ implies $s_\gamma < 2 \max\{\lfloor s_\gamma \rfloor, 1\}$, it remains to show $W_{\mathscr B} \geq \prod_{\gamma \in \mathscr B} \max\{\lfloor s_\gamma \rfloor, 1\}$.
    Let $\gamma \in \mathscr B$. 
    Then, since $W_{\mathscr B} \geq 1$, we must have that
    \begin{equation}
        \left\lfloor \frac{\gamma_+^\top u - \gamma_-^\top \l}{T} \right\rfloor - \left\lceil \frac{\gamma_+^\top \l - \gamma_-^\top u}{T} \right\rceil + 1 \geq 1.
    \end{equation}
    Moreover, using the properties of the ceiling and floor functions
    \begin{equation}
        \left\lfloor \frac{\gamma_+^\top u - \gamma_-^\top \l}{T} \right\rfloor - \left\lceil \frac{\gamma_+^\top \l - \gamma_-^\top u}{T} \right\rceil + 1 > \frac{\gamma_+^\top u - \gamma_-^\top \l}{T} - \frac{\gamma_+^\top \l - \gamma_-^\top u}{T} - 1= s_\gamma - 1,
    \end{equation}
    and we conclude since the left-hand side is an integer.
\end{proof}

The hypothesis $W_{\mathscr B} \geq 1$ is satisfied for all feasible PESP instances: If $W_{\mathscr B} = 0$, then $R$ and hence $Z$ do not contain an integer point. 
If the bounds $\l$ and $u$ are integral, having assumed $u - \l > \mathbf{0}$, then $\varepsilon \geq 1$. 
However, there is also a large class of practical PESP instances, where each arc $a \in A(G)$ is either fixed $(\l_a = u_a)$ or free $(u_a - \l_a = T-1)$, coined \emph{reduced} instances in \cite{patzold_matching_2016}. 
As the former can be contracted, such that only the free arcs remain, while the cyclomatic number $\mu$ is unchanged, \Cref{thm:width-approx} then yields that the minimum width of an integral cycle basis is at least $|\mathscr S(G)|\left( \frac{T-1}{2T} \right)^\mu$ for these instances.
We finally want to remark that the number of spanning trees and the minimum width of an integral cycle basis are exponential in the input size of a PESP instance, so that it makes sense to compare their logarithms. 
Then \Cref{thm:width-approx} states that the logarithms of both numbers differ by $\mu \log(\varepsilon/(2T))$.

\section{Outlook}
\label{sec:outlook}

In order to assess the practical use of tropical methods described in \Cref{sec:tropical} for the optimization of periodic timetables, the next task is to implement and to evaluate the {\tns} heuristic.

A challenging question is how the insights on the cycle offset zonotope (\Cref{sec:zonotopes}) can be exploited for the purpose of optimization. For example, we expect that these zonotopes are related to the space of feasible solutions of a reformulation of PESP by means of Benders' decomposition \cite{LindnervanLieshout2021}. The estimates on the cycle basis width in \Cref{subsec:width} motivate further investigation of approximation algorithms or hardness results for the minimum width cycle basis problem.

Finally, the understanding of cycle offset zonotopes may be improved. It would be interesting to characterize all fine zonotopal tilings in terms of the polytropal decomposition, or by constructions similar to the one in \Cref{subsec:construction}.

\section*{Acknowledgements}

We want to thank Ralf Borndörfer for encouraging us to work on the link between periodic timetabling and tropical geometry.

Enrico Bortoletto has been funded within the Research Campus MODAL, funded by the German Federal Ministry of Education and Research (BMBF) (fund number 05M20ZBM).

Berenike Masing has been funded by the Deutsche Forschungsgemeinschaft (DFG, German Research Foundation) under Germany's Excellence Strategy – The Berlin Mathematics Research Center MATH+ (EXC-2046/1, project ID: 390685689).

\newpage

\bibliographystyle{acm}
\bibliography{references}

\end{document}